\documentclass[reqno]{amsart}
\usepackage{amssymb}
\usepackage{amsthm}
\usepackage{amsmath}
\usepackage{latexsym}
\usepackage[v2,tips]{xy}
\usepackage{mathrsfs}
\usepackage[T1]{fontenc}
\usepackage[latin1]{inputenc}
\usepackage{enumitem}
\usepackage[a4paper, total={5.7in, 9in}]{geometry}
\usepackage{graphicx}
\usepackage{hyperref}
\usepackage[normalem]{ulem}
\usepackage{dutchcal}
\usepackage{tikz-cd} 
\usepackage{cleveref}
\usepackage{tabularx}

\newcommand{\N}{\mathbb{N}}

\newcommand{\R}{\mathbb{R}}

\newcommand{\K}{\mathbb{K}}

\newcommand{\Ss}{\mathcal{S}_{1}}

\newtheorem{theorem}{Theorem}[section]
\newtheorem{lemma}[theorem]{Lemma}
\newtheorem{proposition}[theorem]{Proposition}
\newtheorem{corollary}[theorem]{Corollary}
\theoremstyle{definition}
\newtheorem{definition}[theorem]{Definition}

\newtheorem{remark}[theorem]{Remark}
\newtheorem{construction}[theorem]{Construction}
\newtheorem*{organization}{Organization and overview of content}

\numberwithin{equation}{section}

\renewcommand{\ge}{\geqslant}
\renewcommand{\geq}{\geqslant}
\renewcommand{\le}{\leqslant}
\renewcommand{\leq}{\leqslant}

\newcommand{\gl}{\ensuremath{\Gamma}\!\operatorname{L}_1}
\newcommand{\SC}{\operatorname{SC}}
\newcommand{\lv}{\lVert}
\newcommand{\rv}{\rVert}

\title[Operator ideals on the Baernstein and Schreier spaces]{Closed ideals of operators on the Baernstein\\ and Schreier spaces}
\author[N.J.~Laustsen]{Niels Jakob Laustsen}
\address{(N.J.~Laustsen) School of Mathematical Sciences, Fylde
  College, Lancaster University, Lancaster LA1 4YF, United Kingdom}
\email{n.laustsen@lancaster.ac.uk}

\author[J.~Smith]{James Smith}
\address{(J.~Smith) School of Mathematical Sciences, Fylde
  College, Lancaster University, Lancaster LA1 4YF, United Kingdom}
\email{j.smith43@lancaster.ac.uk}


\begin{document}

\begin{abstract}
We study the lattice of closed ideals of bounded operators on two families of Banach spaces: the Baernstein spaces~$B_p$ for $1<p<\infty$ and  the Schreier spaces~$S_p$ for $1\le p<\infty$. Our main conclusion is that there are~$2^{\mathfrak{c}}$ many closed ideals that lie between the ideals of compact and  strictly singular operators on each of these spaces, and also~$2^{\mathfrak{c}}$ many closed ideals that contain projections of in\-finite rank. 

Counterparts of results of Gasparis and Leung using a numerical index to distinguish the iso\-mor\-phism types of subspaces spanned by subsequences of the unit vector basis for the higher-order Schreier spaces play a key role in the proofs, as does the John\-son--Schecht\-man technique for constructing~$2^{\mathfrak{c}}$ many closed ideals of operators on a Banach space. 
\end{abstract}

\maketitle

\section{Introduction}

\noindent Recently, the lattice of closed ideals of the Banach algebra~$\mathscr{B}(X)$ of bounded operators on a Banach space~$X$ has been studied intensively, in many cases leading to the conclusion that it has cardinality~$2^{\mathfrak{c}}$, which is the largest possible value 
when~$X$ is separable, and an order structure that is at least as complex as the power set of~$\R$. 
We add two families of Banach spaces to the list for which these conclusions can be drawn: the Baernstein spaces~$B_p$ for $1<p<\infty$ and the Schreier spaces~$S_p$ for $1\le p<\infty$. 

Their definitions rely on the family~$\mathcal{S}_1$ of \emph{Schreier sets,} that is, the finite subsets of the natural numbers whose minimum dominates their cardinality. 

We can now define the \emph{$p^{\text{th}}$ Schreier space}~$S_p$, for $1\le p<\infty$, as the completion of the vector space~$c_{00}$  of finitely supported elements $x=(x(n))_{n\in\N}\in\K^\N$ with respect to the norm 
\begin{equation}\label{Eq:Schreiernorm} \lVert x\rVert_{S_p} = \sup\biggl\{\Bigl(\sum_{n\in F} \lvert x(n)\rvert^p\Bigr)^{\frac1p} : F\in\mathcal{S}_1\setminus\{\emptyset\}\biggr\}, \end{equation}
while the \emph{$p^{\text{th}}$ Baernstein space}~$B_p$, for $1<p<\infty$, is the completion of~$c_{00}$ with respect to the norm 
\begin{multline}\label{Eq:Baernsteinnorm} \lVert x\rVert_{B_p} = \sup\biggl\{\biggl(\sum_{j=1}^m\Bigl(\sum_{n\in F_j} \lvert x(n)\rvert\Bigr)^p\biggr)^{\frac1p} : m\in\N,\ F_1,\ldots,F_m\in\mathcal{S}_1\setminus\{\emptyset\}\ \text{and}\biggr.\\[-4mm] 
\biggl.\max F_j<\min F_{j+1}\ \text{for}\ 1\le j<m\biggr\}. 
\end{multline}
(Note that~\eqref{Eq:Baernsteinnorm} would simply define the $\ell_1$-norm for $p=1$, which is why there is no Baernstein space~$B_1$.)

Having introduced these Banach spaces, let us summarize our main conclusions about them. 

\begin{theorem}\label{thmsmalllargeideals}
    Let $E=B_{p}$ for some $1<p<\infty$ or $E=S_{p}$ for some $1\le p <\infty$. Then:
    \begin{enumerate}[label={\normalfont{(\roman*)}}]
        \item\label{itemsmall} $\mathscr{B}(E)$ contains $2^{\mathfrak{c}}$ many closed ideals between the ideals of compact and strictly singular operators. 
    \item\label{itemcompactss} The ideals of strictly singular and inessential operators on $E$ are equal.\\    
    For $E=S_p$, these ideals are also equal to the ideal of weakly compact operators.
    \item\label{itemlarge} $\mathscr{B}(E)$ contains $2^{\mathfrak{c}}$ many closed ideals which are larger than the ideal 
    \[ \{ UV : U\in\mathscr{B}(D,E),\, V\in\mathscr{B}(E,D) \} \] 
    of operators factoring through~$D$, where $D=\ell_p$ if $E=B_p$ and $D=c_0$ if $E=S_p$. 
    \item\label{itemmaxideals} $\mathscr{B}(E)$ contains at least continuum many maximal ideals. 
\end{enumerate}
\end{theorem}

The Baernstein and Schreier spaces originate in Schreier's counter\-example~\cite{schreier} from~1930, which showed that $C[0,1]$ does not have the Banach--Saks property. More precisely, Schreier defined what we call a Schreier set, but did not explicitly consider any of the Banach spaces we study. More than 40 years later,  Baernstein~\cite{AB} introduced the space~$B_2$ to provide an example of a reflexive Banach space without the Banach--Saks property. Not long after, Seifert observed in his dissertation~\cite{S} that Baernstein's definition carries over to arbitrary \mbox{$p\in(1,\infty)$}. 

Surprisingly, it appears that the Schreier space~$S_1$ was not defined until seven years after Baernstein's work~\cite{AB}. It is hard to imagine that Baernstein did not know~$S_1$, but as far as we have been able to find out, Beauzamy~\cite{beau} was the first person to define it explicitly, using it in combination with interpolation methods to obtain another example of a reflexive Banach space without the Banach--Saks property. Bird and the first author~\cite{BL} studied the spaces~$S_p$ for $p>1$. 

A much more commonly researched variant of the Banach space~$S_1$ is the family of \emph{higher-order Schreier spaces}~$X[\mathcal{S}_\xi]$, defined by Alspach and Argyros~\cite{AA} for every countable ordinal~$\xi$; the correspondence is that the space we denote~$S_1$ is equal to~$X[\mathcal{S}_1]$ (and $X[\mathcal{S}_0] = c_0$). We shall not add to the theory of these spaces for $\xi\ge 2$, but instead develop counterparts of some of the main results about them for the Baernstein and Schreier spaces, as we shall explain next.

\begin{organization} We begin by collecting some preliminary material in \Cref{S:2}, most importantly a quantitative version of the fact that for every $1<p<\infty$, the Baern\-stein space~$B_p$ is saturated with complemented copies of~$\ell_p$, while the Schreier spaces are saturated with complemented copies of~$c_0$.
\Cref{thmsmalllargeideals}\ref{itemcompactss} follows easily from these results, as we shall show in \Cref{S:3}. 

The remaining parts of \Cref{thmsmalllargeideals} are substantially harder to verify. \Cref{S:4} contains the main technical tool that we require: a counterpart for the Baernstein and Schreier spaces of some results of Gasparis and Leung~\cite{GL} concerning the higher-order Schreier spaces. For every $n\in\N$, they introduced a numerical index which characterizes when two subspaces spanned by infinite subsequences of the unit vector basis for~$X[\mathcal{S}_n]$ are isomorphic. Surprisingly, we find that their  index for $n=1$ works for the Baernstein spaces~$B_p$ and --- perhaps less surprisingly --- the Schreier spaces~$S_p$ for $p>1$. 

Beanland, Kania and the first author~\cite{BKL} used the results from~\cite{GL} to demonstrate that the family of closed ideals of~$\mathscr{B}(X[\mathcal{S}_n])$ that are singly generated by basis projections has a very rich  structure for every $n\in\N$. In \Cref{S:5}, we show that by referring to \Cref{S:4} instead of~\cite{GL}, we can transfer the arguments from~\cite{BKL} to the Baernstein and Schreier spaces; \Cref{manychainsofideals} states our main conclusions, which include \Cref{thmsmalllargeideals}\ref{itemmaxideals}.

Answering a question raised in~\cite{BKL}, Manoussakis and Pelczar-Barwacz~\cite{MPB} combined the results from~\cite{GL} with the seminal idea of Johnson and Schechtman~\cite{JS} to prove that~$\mathscr{B}(X[\mathcal{S}_n])$ contains~$2^{\mathfrak{c}}$ many closed ideals that lie between the ideals of compact and strictly singular operators  for every $n\in\N$. \Cref{thmsmalllargeideals}\ref{itemsmall} is the analogue of this result for the Baernstein and Schreier spaces. We prove it in \Cref{S:6}, together with \Cref{thmsmalllargeideals}\ref{itemlarge}, whose proof turns out to be the easier of the two. The reason is that \Cref{thmsmalllargeideals}\ref{itemsmall} requires the non-tri\-vial fact that the formal inclusion map from~$B_p$ into~$\ell_p$ is strictly singular. We verify it using an in\-equality due to Jameson, who has generously allowed us to include his proof of it in \Cref{App:GJOJ}.
In contrast to Manoussakis and Pelczar-Barwacz, we express our arguments in terms of the numerical index of Gasparis and Leung, thereby elucidating their combinatorial nature and providing a blueprint for other Banach spaces admitting a suitable index.  
\end{organization}

To provide additional context and background for our results, we conclude this introduction with a survey of separable Banach spaces~$X$ for which the Banach algebra~$\mathscr{B}(X)$ contains~$2^{\mathfrak{c}}$ many closed ideals. As far as we know, Gowers' hyperplane space~$X_{\text{G}}$ originally introduced in~\cite{Gow} is the first example of this kind; more precisely, the first author \cite[Theorem~8.4]{Lau} classified the maximal ideals of~$\mathscr{B}(X_{\text{G}})$ and noted that there are~$2^{\mathfrak{c}}$ of them. 

A major breakthrough occurred when Johnson and Schechtman~\cite{JS} showed that~$\mathscr{B}(L_p[0,1])$ contains~$2^{\mathfrak{c}}$ many closed ideals for every $p\in(1,2)\cup(2,\infty)$. Their key technique has proved very versatile and spawned many new results. \Cref{thm:FSZ} states a variant of it, formulated by Freeman, Schlumprecht and Zs\'{a}k~\cite{FSZ}, who used it to verify that~$\mathscr{B}(X)$ contains~$2^{\mathfrak{c}}$ many closed ideals for a number of direct sums of Banach spaces, notably $X = \ell_p\oplus\ell_q$, $X= \ell_q\oplus c_0$ and $X = \ell_q\oplus\ell_\infty$ for $1\le p<q<\infty$, as well as the Hardy space~$H_1$ and its pre\-dual~$\operatorname{VMO}$. 

Also building on the Johnson--Schechtman technique, Manoussakis and Pelczar-Barwacz~\cite{MPB} showed that $\mathscr{B}(X)$ contains~$2^{\mathfrak{c}}$ many closed ideals for Schlumprecht's arbitrarily distortable Banach space~\cite{Sch} and the higher-order Schreier spaces~$X[\mathcal{S}_n]$  for $n\in\N$, as already mentioned. In collaboration with Causey, Pelczar-Barwacz~\cite{CPB} has subsequently extended the latter result to the Schreier spaces~$X[\mathcal{S}_\xi]$ of any countable order~$\xi$, as well as their duals and biduals. 

Finally, Chu and Schlumprecht~\cite{ChuSch} have shown that~$\mathscr{B}(T[\mathcal{S}_\xi,\theta])$ contains~$2^{\mathfrak{c}}$ many closed ideals for every countable ordinal~$\xi$ and $0<\theta<1$, where $T[\mathcal{S}_\xi,\theta]$ denotes the Tsirelson space of order~$\xi$, as defined by Alspach and Argyros~\cite{AA}.

\section{Preliminaries, including a saturation result for the Baernstein and Schreier spaces}\label{S:2}

\noindent We begin with some general conventions. All vector spaces are over the same scalar field~$\mathbb{K}$, either the real or the complex numbers. We use the letters $X,Y,\ldots$ to denote generic Banach spaces, while we reserve the letter~$E$ for either the Baernstein space~$B_p$ or the Schreier space~$S_p$ and the letter~$D$ for either~$\ell_p$ or~$c_0$, in the same way as in the statement of \Cref{thmsmalllargeideals}. In line with these conventions, $(e_n)_{n\in\N}$ will always denote the unit vector basis for~$E$ (to be discussed in more detail below) and~$(d_n)_{n\in\N}$ the unit vector basis for~$D$. 

The term ``operator'' means a bounded, linear map between two Banach 
spaces~$X$ and~$Y$. We write $\mathscr{B}(X,Y)$ for the space of operators $X \rightarrow Y$ and abbreviate  $\mathscr{B}(X,X)$ to~$\mathscr{B}(X)$ in line with standard practice. 

Let $(x_{n})_{n\in\N}$ be a (Schauder) basis for a Banach space $X$.   For every $n\in\N$, we denote the $n^{\text{th}}$ coordinate functional by~$x_{n}^{*}\in X^*$. Suppose that the basis $(x_{n})$ is un\-con\-di\-tional. Then, for every subset~$N$ of~$\N$, $P_{N}\in\mathscr{B}(X)$ denotes the basis projection given by $P_Nx=\sum_{n\in N}\langle x, x_n^*\rangle x_n$ for $x\in X$. As usual, we abbreviate $P_{\{1,2,\ldots,n\}}$ to~$P_n$ for $n\in\N$.

We follow the convention that $\N=\{1,2,3,\ldots\}$ and write $[\N]$ and $[\N]^{<\infty}$ for the families of infinite and finite subsets of~$\N$, respectively.   As already mentioned, 
\[  \mathcal{S}_{1} = \{F\in[\N]^{<\infty} : \lvert F\rvert\leq \min F\} \]
is the family of \emph{Schreier sets}, where $\lvert F\rvert$ denotes the cardinality of the set~$F$ (and by convention $\min\emptyset = 0$, so $\emptyset\in\mathcal{S}_1$). 
Observe that $\mathcal{S}_{1}$ is closed under taking subsets, and is \emph{spreading} in the following sense: suppose that  $\{ j_{1}<j_2<\cdots<j_{n}\}\in\mathcal{S}_{1}$ and $j_i\le k_i$ for each \mbox{$1\le i\le n$}; then $\{k_{1},\ldots,k_{n}\}\in \mathcal{S}_{1}$. A Schreier set is \emph{maximal} if it is not contained in any strictly larger Schreier set. Clearly, this is equivalent to $F\ne\emptyset$ and $\lvert F\rvert =\min F$.

It will be convenient to express the $p^{\text{th}}$ Schreier and Baernstein norms defined by~\eqref{Eq:Schreiernorm} and~\eqref{Eq:Baernsteinnorm} as the suprema over certain families of semi-norms. This is straightforward for the Schreier norms and was used extensively in~\cite{BL}:  for $1\le p<\infty$ and $x=(x(n))_{n\in\N}\in\K^\N$, we 
can write
\[ \lv x\rv_{S_{p}} = \sup\bigl\{ \mu_{p}(x,F) : F \in \mathcal{S}_{1}\bigr\}\in[0,\infty],\quad \text{where}\quad \mu_p(x,F) = \begin{cases} 0 &\text{if}\ F=\emptyset\\ {\displaystyle{\Bigl(\sum_{n \in F}\lvert x(n)\rvert^{p}\Bigr)^{\frac1p}}} &\text{otherwise.}\end{cases} \] 
As observed in \cite[Lemma~3.3(iv)]{BL}, $\lVert\,\cdot\,\rVert_{S_p}$ defines a complete norm on the subspace $Z_p =\{ x\in\K^\N : \lv x \rv_{S_{p}} < \infty\}$  of~$\K^\N$, and we can view~$S_p$ as the closed subspace of~$Z_p$ spanned by the  
``unit vector basis'' $(e_{n})_{n\in\N}$  given by $e_n(m) = 1$ if $m=n$ and  $e_n(m) = 0$ otherwise. Justifying its name, $(e_{n})_{n\in\N}$ is a 
normalized basis for~$S_p$ 
that is  $1$\nobreakdash-un\-con\-di\-tional and shrinking, as shown in \cite[Prop\-o\-si\-tions~3.5 and~3.10 and Corollary~3.12]{BL}. It is not hard to verify that~$S_p$ is a proper subspace of~$Z_p$; in fact, $Z_p$~is non-separable by \cite[Corollary~5.6]{BL}. 

To analogously express the Baernstein norm as the supremum of a certain family of seminorms, we introduce the following notion: a \emph{Schreier chain} is a non-empty, finite collection~$\mathcal{C}$ of non-empty, consecutive Schreier sets; that is, $\mathcal{C} = \{F_{1}<F_2<\cdots<F_{m}\}$, where $m\in\N$, $F_1,\ldots,F_m\in\Ss\setminus\{\emptyset\}$, and the notation
$F_1<F_2<\cdots<F_m$ signifies that $\max F_{j} < \min F_{j+1}$ for each $1\le j <m$. We write~$\SC$ for the collection of all Schreier chains. Then, for $1<p<\infty$ and a Schreier chain $\mathcal{C}$, we can define a seminorm~$\beta_{p}(\,\cdot\,,\mathcal{C})$ on $\mathbb{K}^{\N}$ by 
\[ \beta_{p}(x,\mathcal{C}) = \biggl(\sum_{F\in \mathcal{C}}\Bigl(\sum_{n\in F}\lvert x(n)\rvert\Bigr)^{p}\biggr)^{\frac{1}{p}}\qquad (x=(x(n))_{n\in\N}\in\K^\N), \]
and 
\[ \lVert x\rVert_{B_{p}} = \sup\bigl\{\beta_{p}(x,\mathcal{C}) :\mathcal{C} \in \SC\bigr\}\in[0,\infty]\qquad (x\in\K^\N). \]
In contrast to the Schreier spaces, it turns out that the $p^{\text{th}}$ Baernstein space (defined as the completion of~$c_{00}$ with respect to this norm) is precisely the collection of vectors $x \in \mathbb{K}^{\N}$ with finite Baernstein norm:  
$B_{p} = \{ x\in\K^\N : \lv x \rv_{B_{p}} < \infty\}$; this follows by replacing the exponent~$2$ with~$p$ in Baernstein's argument given in \cite[page~92, first paragraph of ``Proof of~(2)'']{AB}. As is the case for the Schreier spaces, the unit vector basis $(e_n)_{n\in\N}$ forms a $1$\nobreakdash-un\-con\-di\-tional, 
normalized basis for~$B_p$. 

In line with standard practice, $\operatorname{supp}(x) = \{ n\in\N : x(n)\ne0\}$ denotes the \emph{support} of an element $x=(x(n))_{n\in\N}\in\K^\N$. Clearly, when computing $\lv x\rv_{S_p}$, it suffices to consider~$\mu_p(x,F)$ for $F\in\mathcal{S}_1$ with $F\subseteq\operatorname{supp}(x)$. Similarly, when computing $\lv x\rv_{B_p}$, it suffices to consider $\beta_p(x,\mathcal{C})$ for $\mathcal{C}\in\SC$ with $\bigcup\mathcal{C}\subseteq\operatorname{supp}(x)$.

\begin{lemma}\label{supportcovered}
Let $1<p<\infty$, and suppose that the non-zero coordinates of $x\in c_{00}$ are decreasing in absolute value. Then~$x$ attains its Baernstein norm at some Schreier chain covering $\operatorname{supp}(x);$ that is, $\lv x \rv _{B_{p}} = \beta_p(x,\mathcal{C})$ for some $\mathcal{C}\in\SC$ with $\bigcup\mathcal{C} = \operatorname{supp}(x)$.
\end{lemma}

\begin{proof}
Take $\mathcal{C} = \{F_1<F_2<\cdots<F_m\}\in\SC$ with $\bigcup\mathcal{C}\subsetneq\operatorname{supp}(x)$. For $1 \leq j \leq m$, let $F_{j}'$ be the set which contains precisely the first $\lvert F_{j}\rvert$ points of $\operatorname{supp}(x)\cap [\min(F_{j}),\max(F_{j})]$. Then $\mathcal{C}' = \{F_1'<F_2'<\cdots<F_m'\}\in\SC$ and
\begin{equation*}
    \sum_{n\in F_{j}}\lvert\langle x, e_{n}^*\rangle\rvert\leq \sum_{n\in F_{j}'}\lvert\langle x, e_{n}^*\rangle\rvert\qquad (1 \leq j \leq m)
\end{equation*}
because the non-zero coordinates of~$x$ are decreasing in absolute value. Hence $\beta_p(x,\mathcal{C})\le\beta_p(x,\mathcal{C}')$. 
The set 
$\mathcal{C}'' = \bigl\{ \{n\} : n\in\operatorname{supp}(x) \setminus\bigcup\mathcal{C}'\bigr\}$ is non-empty because $\bigcup\mathcal{C}\subsetneq\operatorname{supp}(x)$. Clearly $\mathcal{D} = \mathcal{C}'\cup\mathcal{C}''$ is a Schreier chain, and 
\[ \beta_{p}(x,\mathcal{C})\le \beta_{p}(x,\mathcal{C}') < \beta_{p}(x,\mathcal{D})\le\lv x\rv_{B_p}, \] 
so $\lv x\rv_{B_p}$ is not attained at~$\mathcal{C}$. However, it must be attained at some Schreier chain because~$x$ is finitely supported.
\end{proof}

For later reference, we now record an estimate for the norm of certain ``flat'' vectors. 

\begin{lemma}\label{flatvectorestimate}
Let $\{F_1<F_2<\cdots<F_m\}$ be a chain of maximal Schreier sets for some $m\in\N$.  Then
\begin{alignat}{3}
    1&\leq \biggl\lVert \sum_{j=1}^{m}\frac{1}{|F_{j}|^{\frac{1}{p}}}\sum_{k\in F_{j}}e_{k}\biggr\rVert_{S_{p}} &&\leq 2^\frac{1}{p}&\qquad &(1\le p<\infty)\label{estimate2schreier}\\
\intertext{and}
m^{\frac{1}{p}}&\le \biggl\lVert \sum_{j=1}^{m}\frac{1}{|F_{j}|}\sum_{k\in F_{j}}e_{k}\biggr\rVert_{B_p} &&\le 2m^{\frac{1}{p}}&\qquad &(1< p<\infty).\label{estimate2baernstein} 
\end{alignat}
\end{lemma}
\begin{proof} To prove~\eqref{estimate2schreier}, set $x = \sum_{j=1}^{m}\lvert F_{j}\rvert^{-\frac{1}{p}}\sum_{k \in F_{j}}e_{k}\in S_p$. By definition, we have $\mu_p(x,F_j) = 1$ for $1\le j\le m$; the lower bound on the norm follows. On the other hand, since~$x$ is finitely supported and its non-zero coordinates are decreasing, $\lv x \rv_{S_{p}}$ is attained at some set $G\in\mathcal{S}_1$ that intersects at most two consecutive sets from $\{F_{1},\ldots ,F_{m}\}$; that is, $G\subseteq F_j\cup F_{j+1}$ for some $1\le j<m$, and we have $\mu_{p}(x,G)^{p}\leq\mu_p(x,F_{j})^p + \mu_p(x,F_{j+1})^p = 2$, which proves the upper bound. 

Turning our attention to~\eqref{estimate2baernstein}, we define $x = \sum_{j=1}^{m}\frac{1}{|F_{j}|}\sum_{k\in F_{j}}e_{k}\in B_p$. The lower bound on the norm of~$x$ follows from the fact that  $\mathcal{C} = \{F_1<F_2<\cdots<F_m\}$ is a Schreier chain for which $\beta_{p}(x,\mathcal{C}) = m^{\frac{1}{p}}$. 

For the upper bound, \Cref{supportcovered} implies that we can take $\mathcal{C}=\{G_1<\cdots<G_n\}\in\SC$ such that $\lVert x\rVert_{B_p} = \beta_p(x,\mathcal{C})$ and $\bigcup_{k=1}^{n}G_{k} = \bigcup_{j=1}^{m}F_{j}$. Define a map $\varphi\colon \{1,\ldots,n\}\to\{1,\ldots,m\}$ by 
\[ \varphi(k) = \min\{j : G_{k} \cap F_{j} \neq \emptyset\}. \]
We observe that~$\varphi$ is surjective because otherwise we would have $F_j\subsetneq G_k$ for some~$j$ and~$k$, contradicting that~$F_j$ is a maximal Schreier set. 

Fix $j\in\{1,\ldots,m\}$, set $h(j)=\max\varphi^{-1}(\{j\})$, and note that $\bigcup\bigl\{ G_k : \varphi(k)=j,\, k\ne h(j)\bigr\}\subseteq F_j$ because the sets $G_1,\ldots,G_n$ are successive. Since the $\ell_1$-norm dominates the $\ell_p$-norm, it follows that 
\[ \sum_{k\in\varphi^{-1}(\{j\})}\Bigl(\sum_{i\in G_k}\lvert \langle x, e_i^*\rangle\rvert\Bigr)^p\le  \biggl(\sum_{k\in\varphi^{-1}(\{j\})}\sum_{i\in G_k}\lvert \langle x, e_i^*\rangle\rvert\biggr)^p\le \biggl(1+\frac{\lvert G_{h(j)}\cap F_{j+1}\rvert}{\lvert F_{j+1}\rvert}\biggr)^p\le 2^p. \] 
Combining this with the fact that the sets $\{\varphi^{-1}(\{j\}) : 1\le j\le m\}$ partition~$\{1,\ldots,n\}$, we conclude that
\[ \lVert x\rVert_{B_p}^p = \beta_p(x,\mathcal{C})^p = \sum_{j=1}^m\sum_{k\in\varphi^{-1}(\{j\})}\Bigl(\sum_{i\in G_k}\lvert\langle x,e_i^*\rangle\rvert\Bigr)^p\le 2^pm, \]
which gives the upper bound. 
\end{proof}

We shall now present the main conclusion of this section, which is a quantitative version of two results in the literature: one by Seifert \cite[Theorem~3]{S} stating that~$B_p$ is saturated with complemented copies of~$\ell_{p}$ for \mbox{$1<p<\infty$}, the other by Bird and the first author \cite[Corollary~5.4]{BL} stating that the Schreier spaces are saturated with copies of~$c_{0}$ (which are automatically complemented by Sobczyk's Theorem). Our statement strengthens these results by providing explicit norm bounds on the projections and isomorphisms, using the following terminology.

\begin{definition}\label{defnCsaturated}
    Let~$X$ and~$Y$ be Banach spaces. We say that~$X$ is \textit{$C$-uniformly saturated with complemented copies of~$Y$} for some constant $C \geq 1$ if every closed, infinite-dimensional sub\-space of~$X$ con\-tains a closed sub\-space~$Z$ for which
    \begin{enumerate}[label={\normalfont{(\roman*)}}]
        \item\label{defnCsaturated:i} there exists an isomorphism~$U$ of~$Y$ onto~$Z$ with $\lv U \rv\cdot \lv U^{-1} \rv \leq C$, and
        \item\label{defnCsaturated:ii} there exists a projection~$P$ of~$X$ onto~$Z$ with $\lv P \rv \leq C$.
    \end{enumerate}
\end{definition}

\begin{theorem}\label{lpsaturation}
    Let $(E,D) = (B_{p},\ell_{p})$ for some $1 < p < \infty$ or $(E,D) = (S_{p},c_{0})$ for some $1 \leq p < \infty$. Then~$E$ is $C$-uniformly saturated with complemented copies of $D$ for every $C > 1$.
\end{theorem}

We provide a full proof of this theorem, even though the norm bounds it provides may seem only a modest improvement of existing knowledge. However, Seifert never published his dissertation, and it contains an unfortunate error: \cite[Lemma~2]{S} claims that  every semi-normalized block basic sequence $(w_n)_{n\in\N}$ of the unit vector basis for~$B_p$ admits a subsequence which is equivalent to the unit vector basis for~$\ell_p$. That is not true; for instance, no subsequence of the unit vector basis for~$B_p$ is equivalent to the $\ell_p$-basis. The correct statement is that $(w_n)_{n\in\N}$ admits a block basic sequence which is equivalent to the $\ell_p$-basis; this will follow from \Cref{infinitynull} and \Cref{P:charDsubseq} below.

To compound this issue, Seifert's incorrect statement has been reproduced in  \cite[Theorem~0.15(c)--(d)]{CS}, as well as \cite[page~233]{fgdb} in a special case. This has caused at least one mis\-take in the published literature: citing~\cite{CS}, Flores \emph{et al.}\ \cite[page~334]{FHST} deduce that the Baern\-stein spaces are ``dis\-joint\-ly homo\-geneous''. However, that is impossible because by \cite[Theorem~2.13]{FHST}, it would imply that every strictly singular operator on~$B_p$ is compact, which would contradict \Cref{thmsmalllargeideals}\ref{itemsmall}. (In more concrete terms, it would also contradict \Cref{P:formalinclusionmap} below, which implies that the formal inclusion map of $B_p$ into~$\ell_p$ composed with any isomorphic embedding of~$\ell_p$ into~$B_p$ is a strictly singular, non-compact operator on~$B_p$.)

Having explained \emph{why} \Cref{lpsaturation} requires a detailed proof, we shall now present one, proceeding through a series of lemmas. 

\begin{definition} A \emph{block subspace} of a Banach space~$X$ with a basis $(x_{n})_{n\in\N}$ is  a closed subspace of~$X$ of the form $\overline{\operatorname{span}}\,(w_{n}:n \in \N)$ for some block basic sequence $(w_{n})_{n\in\N}$ of $(x_{n})_{n\in\N}$. 
\end{definition}

\begin{lemma}\label{C2afterC1}
    Let $X$ and $Y$ be Banach spaces, where $X$ has a basis $(x_{n})_{n\in\N}$. Suppose that there is a constant $C_{1} \geq 1$ for which every block subspace $W$ of $X$ admits operators $U\in\mathscr{B}(Y,W)$ and $V\in\mathscr{B}(X,Y)$ such that $V|_{W}U = I_{Y}$ and $\lv U \rv \cdot \lv V \rv \leq C_{1}$. Then~$X$ is $C_{2}$\nobreakdash-uni\-formly saturated with complemented copies of~$Y$ for every constant $C_{2} > C_{1}$.
\end{lemma}
\begin{proof}
    Let $K$ be the basis constant of $(x_{n})_{n\in\N}$, and let $(P_{n})_{n \in \N}$ be the corresponding basis projections. Given $C_{2} > C_{1}$, choose $\varepsilon \in (0,1)$ such that 
    \begin{equation}\label{C2afterC1:eq2} 
    \frac{7(C_{2}-C_{1})}{4(C_{1}+C_{2})} \geq \varepsilon. 
    \end{equation}
    
    Set $m_{0} = 0$ and $P_{0} = 0$, and  let $Z$ be a closed, infinite-dimensional subspace of $X$. By recursion, we can choose natural numbers $m_{1} < m_{2} < \cdots$ and unit vectors $z_{n} \in Z\cap \ker P_{m_{n-1}}$ such that $\lv z_{n} - P_{m_{n}}z_{n} \rv \leq\varepsilon/(2^{n+2}K)$ for every $n \in \N$. Set $w_{n} = P_{m_{n}}z_{n} \in X$, and note that
    \begin{equation}\label{C2afterC1:eq1}
        \lv w_{n} \rv = \lv z_{n} - (z_{n}-w_{n}) \rv \geq 1 - \frac{\varepsilon}{2^{n+2}K} \geq \frac{7}{8}\qquad (n \in \N).
    \end{equation}
    In particular, $w_{n} \neq 0$, and since $z_{n} \in \operatorname{ker}P_{m_{n-1}}$, it follows that $(w_{n})_{n \in \N}$ is a block basic sequence of $(x_{n})_{n\in\N}$.
    Set $W = \overline{\operatorname{span}}\,(w_{n}: n \in \N)$. By hypothesis, we can find operators $U_{1}\in\mathscr{B}(Y,W)$ and $V_{1}\in\mathscr{B}(X,Y)$ such that $V_{1}|_{W}U_{1} = I_{Y}$ and $\lv U_{1} \rv \cdot \lv V_{1} \rv \leq C_{1}$. 
    
    For each $n \in \N$, choose a functional $f_{n} \in X^{*}$ of norm 1 such that $\langle w_{n},f_{n} \rangle = \lv w_{n} \rv$. Then, using~\eqref{C2afterC1:eq1}, we have
    \begin{gather*}
        \sum_{n=1}^{\infty}\frac{\lv (P_{m_{n}}-P_{m_{n-1}})^{*}f_{n} \rv}{\lv w_{n} \rv}\cdot \lv w_{n} - z_{n} \rv \leq \sum_{n=1}^{\infty}\frac{2K}{{7}/{8}}\cdot \frac{\varepsilon}{2^{n+2}K} = \frac{4\varepsilon}{7},
    \end{gather*}
    so we can define an operator $R\in\mathscr{B}(X)$ by
    \begin{gather*}
        R = \sum_{n=1}^{\infty}\frac{(P_{m_{n}} - P_{m_{n-1}})^{*}f_{n}}{\lv w_{n} \rv} \otimes (w_{n}-z_{n}),
    \end{gather*}
    where $f\otimes x$, for $f\in X^*$ and $x\in X$, denotes the rank-one operator $y\mapsto \langle y,f\rangle x$ as usual.
    Since $\lv R \rv\le 4\varepsilon/7 < 1$, the Neumann series implies that the operator $S = I_{X}- R \in \mathscr{B}(X)$ is invertible, and $\lv S^{-1} \rv \leq (1 - 4\varepsilon/7)^{-1} = 7/(7-4\varepsilon)$. The definition of $R$ shows that    
    \begin{gather*}
        Rw_{j} = \sum_{n=1}^{\infty} \frac{\langle (P_{m_{n}} - P_{m_{n-1}})w_{j},f_{n}\rangle}{\lv w_{n} \rv}(w_{n}-z_{n}) = \frac{\langle w_{j},f_{j} \rangle}{\lv w_{j} \rv}(w_{j}-z_{j}) = w_{j}- z_{j}, 
    \end{gather*}
    so $Sw_{j} = z_{j}$ for each $j \in \N$, and therefore $S[W]\subseteq Z$. It follows that $Z_0 = (S|_WU_1)[Y]$ is a subspace of~$Z$, and the operators $U = S|_{W}U_{1}\in \mathscr{B}(Y,Z_0)$ and $V = V_{1}S^{-1}\in \mathscr{B}(X,Y)$ satisfy $V|_{Z_0}U = V_{1}|_{W}U_{1} = I_{Y}$. Since~$U$ is surjective by definition, this implies that~$U$ is invertible with inverse $V|_{Z_0}$ and $P = UV$ is a projection of~$X$ onto~$Z_0$. In particular, $Z_0$ is a closed subspace of~$Z$, and the norm bounds on $\lv U\rv\cdot \lv U^{-1}\rv$ and~$\lv P\rv$ specified in \Cref{defnCsaturated}\ref{defnCsaturated:i}--\ref{defnCsaturated:ii} follow from the fact that  
    \[ \lv U \rv \cdot \lv V \rv \leq \lv S \rv \cdot \lv U_{1} \rv \cdot\lv V_{1} \rv \cdot \lv S^{-1} \rv \leq \Bigl(1 + \frac{4\varepsilon}{7}\Bigr)C_{1} \cdot \frac{7}{7-4\varepsilon} \leq C_{2}, \] where the final inequality is a direct consequence of~\eqref{C2afterC1:eq2}. 
\end{proof}

It turns out that the supremum norm $\lv x\rv_{\infty} = \sup_{n\in\N}\lvert x(n)\rvert$ for $x=(x(n))_{n\in\N}\in\K^\N$ plays an important auxiliary role in a number of results about the Baernstein and Schreier spaces. We shall sometimes use the coordinate functionals to express it in the alternative form $\lv x\rv_{\infty} = \sup_{n\in\N}\lvert\langle x,e_n^*\rangle\rvert$ for $x\in B_p$ or $x\in S_p$.

\begin{lemma}\label{infinitynull}
Every block basic sequence of the unit vector basis for~$B_{p}$ \emph{(}for \mbox{$1 < p < \infty)$} or~$S_{p}$ \emph{(}for $1 \leq p < \infty)$ admits a normalized block basic sequence $(u_{n})_{n\in\N}$ for which $\lv u_{n} \rv_{\infty} \rightarrow 0$ as $n \rightarrow \infty$.
\end{lemma}

\begin{proof}
As usual, let $E=B_p$ or $E=S_p$, and let $(w_n)_{n\in\N}$ be a  block basic sequence of the unit vector basis~$(e_n)_{n\in\N}$ for~$E$. Replacing $(w_n)_{n\in\N}$ with the block basic sequence $(w_n/\lv w_{n}\rv_E)_{n\in\N}$, we may suppose that $(w_n)_{n\in\N}$ is normalized in the $E$-norm.
If $(w_n)_{n\in\N}$ admits a subsequence $(w_{n_j})_{j\in\N}$ such that $\lv w_{n_j} \rv_{\infty} \rightarrow 0$ as $j\to\infty$, there is nothing to prove. 

Otherwise $\delta := \inf_{n\in\N}\lv w_n\rv_{\infty} >0$, so for each $n\in\N$, we can choose $m_n\in\N$ such that $\lvert\langle w_n,e_{m_n}^*\rangle\rvert\ge\delta$. Since $(w_n)_{n\in\N}$ is a block basic sequence, we have $m_1<m_2<\cdots$ and $F_n = \{ m_j :  2^{n-1}\le j< 2^n\}$ is a Schreier set, being a spread of the interval $[2^{n-1},2^n)\cap\N\in\mathcal{S}_1$. This implies that the block basic sequence $(v_n)_{n\in\N}$ of $(w_{n})_{n\in\N}$ defined by  
\[   v_{n} = \sum_{j = 2^{n-1}}^{2^n-1} w_{j}\qquad (n\in\N) \]
is unbounded because 
\[ \lv v_n\rv_{E}\ge \begin{cases} \mu_{p}(v_{n},F_n)\geq 2^{(n-1)/p}\delta\to\infty\quad&\text{as}\quad n\to\infty\quad\text{for}\quad E=S_p,\\
\beta_{p}(v_{n},\{F_n\})\geq 2^{n-1}\delta\to\infty\quad&\text{as}\quad n\to\infty\quad\text{for}\quad E=B_p. \end{cases} \]
On the other hand, $\lv v_{n} \rv_\infty = \max\bigl\{ \lv w_{j}\rv_\infty : 2^{n-1}\le j< 2^n\bigr\}\le 1$ because $(w_j)_{j\in\N}$ is a nor\-mal\-ized block basic sequence of $(e_n)_{n\in\N}$, so $(u_{n} = v_{n}/\lv v_{n} \rv_{E})_{n\in\N}$ is a normalized block basic sequence of~$(w_{n})_{n\in\N}$ such that 
\begin{gather*}
    \lv u_{n} \rv_{\infty} = \frac{\lv v_{n} \rv_\infty}{\lv v_{n} \rv_{E}}\leq \frac{1}{\lv v_{n} \rv_{E}} \rightarrow 0\quad\text{as}\quad n\to\infty. \qedhere
\end{gather*}
\end{proof}

Our next lemma involves the following standard piece of terminology. 
A basic sequence $(x_n)_{n\in\N}$ in a Banach space~$X$ \emph{dominates} a basic sequence $(y_n)_{n\in\N}$ in a Banach space~$Y$ if there is a constant $C>0$ such that
\begin{equation}\label{Eq:Cdomination}
 \biggl\lVert\sum_{n=1}^m \alpha_n y_{n}\biggr\rVert_Y\leqslant C
\biggl\lVert\sum_{n=1}^m \alpha_n x_{n}\biggr\rVert_X\qquad (m\in\N,\,
\alpha_1,\ldots,\alpha_m\in\mathbb{K}). 
\end{equation}
If we wish to record the value of the constant~$C$, we say that~$(x_n)_{n\in\N}$ \emph{$C$-dominates}~$(y_n)_{n\in\N}$.  

\begin{lemma}\label{dominatedbylp}
Let $(E,D) = (B_{p},\ell_{p})$ for some $1 < p < \infty$ or $(E,D) = (S_{p},c_{0})$ for some \mbox{$1 \leq p < \infty$}, and suppose that $(u_{n})_{n\in\N}$ is a normalized block basic sequence of the unit vector basis 
for $E$ with $\inf_{n \in \N}\|u_{n}\|_\infty = 0$. Then, for every constant $C>1$, $(u_n)_{n\in\N}$ admits a subsequence which is $C$-dominated by the unit vector basis 
for~$D$.
\end{lemma}
\begin{proof} Set $\varepsilon = C^{p}-1 > 0$. 
  We begin with the easier case, which is the Schreier space; that is, $(E,D) = (S_{p},c_0)$. 
  Re\-cur\-sive\-ly, we can choose integers $1=j_1<j_2<\cdots$ such that
  \begin{equation}\label{dominatedbylp:eq1}
    \lv u_{j_{k+1}} \rv_{\infty}^p \leq \frac{\varepsilon}{\max(\operatorname{supp}(u_{j_{k}}))}\qquad (k\in\N).
  \end{equation}
  In order to verify that the unit vector basis for~$c_0$ $C$-dominates the basic sequence $(u_{j_k})_{k\in\N}$ in~$S_p$, we  must show that $\mu_p(x,F)\leq C$ whenever $x=\sum_{k=1}^{n}\alpha_ku_{j_{k}}$ for some $n \in \N$ and some $\alpha_1,\ldots,\alpha_n\in\K$ with $\max_{1\le k\le n}\lvert\alpha_k\rvert\le 1$, and $F\in \mathcal{S}_{1}\setminus\{\emptyset\}$ with $F\subseteq \operatorname{supp} (x)$. Set 
  \[ m = \min\{ k\in\N: F\cap\operatorname{supp} (u_{j_k})\ne \emptyset\}. \]
  Then we have $\lvert F\rvert\leq \min F\leq \max(\operatorname{supp} (u_{j_{m}}))$, so $\lvert F\rvert\cdot\lVert u_{j_k}\rVert_{\infty}^p\le\varepsilon$ for $k>m$ by~\eqref{dominatedbylp:eq1}, and therefore
  \begin{align*} 
    \mu_{p}(x,F)^{p} &= \mu_p(\alpha_mu_{j_m},F)^p + \mu_p\biggl(\sum_{k=m+1}^n\alpha_ku_{j_k},F\biggr)^p\\ &\le \lvert\alpha_m\rvert^p\,\lVert u_{j_m}\rVert_{S_p}^p + \lvert F\rvert\Bigl(\max_{m<k\le n}\lvert\alpha_k\rvert\, \lVert u_{j_k}\rVert_{\infty}\Bigr)^p\le 1+\varepsilon=C^p,  
  \end{align*}
  as required.  
    
Proceeding to the case $(E,D) = (B_{p},\ell_p)$, we use the fact that the function $t\mapsto t^p$ is uniformly continuous on~$[0,2]$ to choose numbers $\delta_{k} \in (0,1)$ such that
\begin{equation}\label{dominatedbylp:eq2}
    (s+t)^{p} \leq s^{p} + \frac{\varepsilon}{2^{k}}\qquad (k \in \N,\,s\in[0,1],\,t\in[0,\delta_k]). 
\end{equation}
After replacing $(u_{n})_{n\in\N}$ with a suitable subsequence, we may suppose that $\lv u_{n} \rv_{\infty} \rightarrow 0$ as $n\to\infty$. We can then recursively choose integers $1= j_{1} < j_{2} < \cdots$ such that 
\begin{equation}\label{dominatedbylp:eq3}
    \lv u_{i} \rv_{\infty} \leq \frac{\delta_{k}}{\max(\operatorname{supp}(u_{j_{k}}))}\qquad (k\in\N,\,i \geq j_{k+1}). 
\end{equation}
We seek to verify that the unit vector basis for~$\ell_p$ $C$-dominates the basic sequence $(u_{j_{k}})_{k \in \N}$ in~$B_p$. This amounts to showing that $\beta_p(x,\mathcal{C})\leq C$ whenever $x = \sum_{k=1}^{n}\alpha_{k}u_{j_{k}}$ for some $n \in \N$ and some $\alpha_{1},\ldots,\alpha_{n} \in \mathbb{K}$ with $\sum_{k=1}^{n}\lvert \alpha_{k}\rvert^{p}\leq 1$, and~$\mathcal{C}$ is a Schreier chain contained in $\operatorname{supp}(x)$. Set
\[  \mathcal{C}_{k} = \{F\in \mathcal{C} : \min F\in \operatorname{supp}(u_{j_{k}})\}\qquad (1 \leq k \leq n). \]
Then, defining $\beta_p(x,\emptyset)=0$ to cover the case where $\mathcal{C}_k=\emptyset$ for some~$k$, we can write
\begin{equation}\label{dominatedbylp:eq4}
    \beta_{p}(x,\mathcal{C})^{p} = \sum_{k=1}^{n}\sum_{F\in \mathcal{C}_{k}}\Bigl(\sum_{i\in F}|\langle x,e_{i}^{*}\rangle|\Bigr)^{p} = \sum_{k=1}^{n}\beta_p(x,\mathcal{C}_k)^p.
\end{equation}
We claim that 
\begin{equation}\label{dominatedbylp:eq5}
\beta_p(x,\mathcal{C}_k)^p\le \lvert\alpha_k\rvert^p+\frac{\varepsilon}{2^k}\qquad (1\le k\le n),
\end{equation}
from which the conclusion will follow because substituting~\eqref{dominatedbylp:eq5} into~\eqref{dominatedbylp:eq4}, we obtain   
\[ \beta_{p}(x,\mathcal{C})^{p}\le \sum_{k=1}^n \Bigl(\lvert\alpha_k\rvert^p+\frac{\varepsilon}{2^k}\Bigr)\le 1+\varepsilon =C^p. \]

It remains to prove~\eqref{dominatedbylp:eq5}. Take $k\in\{1,\ldots,n\}$ with $\mathcal{C}_k\ne\emptyset$, let~$G_k$ be the final set in~$\mathcal{C}_k$ (in the sense that~$G_k$ is the set in~$\mathcal{C}_k$ with the largest minimum), and define 
\[ G_{k}' = G_k\cap\operatorname{supp}(u_{j_{k}})\qquad\text{and}\qquad G_{k}'' = G_k\setminus \operatorname{supp}(u_{j_{k}}). \] 
Then we have 
\begin{align}\label{dominatedbylp:eq6} 
\beta_p(x,\mathcal{C}_k)^p &= \sum_{F\in\mathcal{C}_k\setminus\{G_k\}}\!\Bigl(\sum_{i\in F}\lvert\langle x,e_i^*\rangle\rvert\Bigr)^p + \Bigl(\sum_{i\in G_k}\lvert\langle x,e_i^*\rangle\rvert\Bigr)^p\\ &= 
\lvert\alpha_k\rvert^p\!\!\sum_{F\in\mathcal{C}_k\setminus\{G_k\}}\!\Bigl(\sum_{i\in F}\lvert\langle u_{j_k},e_i^*\rangle\rvert\Bigr)^p + (s_k+t_k)^p,\notag 
\end{align}
where we have introduced the quantities
\[ s_k =  \sum_{i\in G_k'}\lvert\langle x,e_i^*\rangle\rvert =  \lvert\alpha_k\rvert\sum_{i\in G_k'}\lvert\langle u_{j_k},e_i^*\rangle\rvert\qquad\text{and}\qquad  t_k = \sum_{i\in G_k''}\lvert\langle x,e_i^*\rangle\rvert. \]
Now we observe that $0\le s_k\le\lvert\alpha_k\rvert\,\lVert u_{j_k}\rVert_{B_p}\le 1$ because $G_k'\in\mathcal{S}_1$, and 
\[ 0\le t_k\le\lvert G_k''\rvert\max_{k<i\le n}\lvert\alpha_i\rvert\,\lVert u_{j_i}\rVert_\infty\le\delta_k, \]
where we have used~\eqref{dominatedbylp:eq3} together with the fact that $\lvert G_{k}''\rvert < \lvert G_k\rvert\leq \min(G_{k})\leq \max(\operatorname{supp}(u_{j_{k}}))$. Hence~\eqref{dominatedbylp:eq2} implies that $(s_k+t_k)^{p} \leq s_k^{p} + \varepsilon/2^k$. Substituting this into~\eqref{dominatedbylp:eq6} and defining the Schreier chain $\mathcal{C}_k' = (\mathcal{C}_k\setminus\{G_k\})\cup\{G_k'\}$, we obtain
\begin{align*} \beta_p(x,\mathcal{C}_k)^p &\le
\lvert\alpha_k\rvert^p\!\!\sum_{F\in\mathcal{C}_k\setminus\{G_k\}}\!\Bigl(\sum_{i\in F}\lvert\langle u_{j_k},e_i^*\rangle\rvert\Bigr)^p +\lvert\alpha_k\rvert^p\Bigl(\sum_{i\in G_k'}\lvert\langle u_{j_k},e_i^*\rangle\rvert\Bigr)^p + \frac{\varepsilon}{2^k}\\ &= \lvert\alpha_k\rvert^p\beta_p(u_{j_k},\mathcal{C}_k')^p + \frac{\varepsilon}{2^k}\le 
\lvert\alpha_k\rvert^p + \frac{\varepsilon}{2^k}. \qedhere
\end{align*}
\end{proof}

\begin{lemma}\label{Sigmaoperator} Let $\mathcal{C} = \{F_{1} < F_{2} < \cdots\}$ be an infinite chain of successive Schreier sets, and take $1<p<\infty$.  Then, for each $x\in B_p$, 
\begin{equation}\label{Sigmaoperator:eq1}
    \Sigma_{\mathcal{C}}x = \Bigl(\sum_{j\in F_{n}}\langle x, e_j^*\rangle\Bigr)_{n\in\N}
\end{equation}
defines an element of~$\ell_p$ with $\lv \Sigma_{\mathcal{C}}x \rv_{\ell_p} \leq \lv x \rv_{B_{p}}$. Hence~\eqref{Sigmaoperator:eq1} defines a map \mbox{$\Sigma_{\mathcal{C}}\colon B_{p} \rightarrow \ell_{p}$}, which is bounded and linear with norm~$1$.
\end{lemma}
\begin{proof}
For $x\in B_p$ and $m\in\N$, we have
\begin{gather*}
\sum_{n=1}^{m}\biggl\lvert\sum_{j\in F_{n}} \langle x,e_{j}^{*}\rangle\biggr\rvert^{p}\leq \sum_{n=1}^{m}\biggl(\sum_{j \in F_{n}}\lvert\langle x,e_{j}^{*}\rangle\rvert\biggr)^{p} = \beta_{p}\bigl(x,\{F_1<F_2<\cdots<F_m\}\bigr)^{p} \leq \lv x \rv_{B_{p}}^p. 
\end{gather*}
This shows that $\Sigma_{\mathcal{C}}x\in\ell_p$ with $\lv \Sigma_{\mathcal{C}}x \rv_{\ell_p} \leq \lv x \rv_{B_{p}}$ because the upper bound~$\lv x \rv_{B_{p}}^p$ is independent of~$m$. The remainder of the lemma is now straightforward to verify.
\end{proof}

\begin{lemma}\label{operatorontolpbasis}  
Let $(E,D) = (B_{p},\ell_{p})$ for some $1 < p < \infty$ or $(E,D) = (S_{p},c_{0})$ for some $1 \leq p < \infty$, and suppose that $(u_{n})_{n\in\N}$ is a normalized block basic sequence  of the unit vector basis for~$E$. Then there exists an operator $V\in\mathscr{B}(E,D)$ of norm~$1$ such that $Vu_{n} = d_{n}$ for every $n \in \N$.
\end{lemma}

\begin{proof}
We begin with the case $(E,D) = (S_{p},c_0)$. For $n \in \N$, set $m_{n} = \max(\operatorname{supp}(u_{n}))$, and  use the Hahn--Banach Theorem to find a functional $f_{n} \in S_p^*$ such that $\langle u_{n},f_{n}\rangle = 1 = \lv f_{n}\rv$. Then, for each $x\in S_{p}$, we can define 
\begin{equation}\label{operatorontolpbasis:eq1}  
    Vx =  \bigl(\langle (P_{m_n}-P_{m_{n-1}})x,f_{n}\rangle\bigr)_{n \in \N}\in\ell_\infty,
\end{equation}  
where we have introduced $m_0=0$ and $P_0=0$ for notational convenience. We see that $Vx\in c_0$ because  
\[ \bigl\lvert\langle (P_{m_n}-P_{m_{n-1}})x,f_{n}\rangle\bigr\rvert\le\lVert (P_{m_n}-P_{m_{n-1}})x\rVert_{S_p}\le \lVert (I-P_{m_{n-1}})x\rVert_{S_p}\to0\quad \text{as}\quad n\to\infty, \]
so~\eqref{operatorontolpbasis:eq1} defines a map $V\colon S_{p} \rightarrow c_0$, which is clearly linear. Furthermore, $V$ is bounded with $\lVert V\rVert\le 1$ because $\lVert P_{m_n}-P_{m_{n-1}}\rVert = 1 = \lVert f_n\rVert$ for $n \in \N$, and we have $Vu_{n} = d_{n}$ because $\operatorname{supp}(u_n)\subseteq (m_{n-1},m_n]$ and $\langle u_{n},f_{n}\rangle = 1$. 

As before, the case $(E,D) = (B_{p},\ell_p)$ is somewhat more involved. We begin by choosing a sequence of scalars $(\sigma_{j})_{j\in\N}$ as follows. If $j \in \operatorname{supp}(u_{n})$ for a (necessarily unique) $n \in \N$, we take $\sigma_{j} \in \mathbb{K}$ of modulus~$1$ such that $\sigma_{j}\cdot\langle u_{n},e_{j}^{*}\rangle > 0$. Otherwise (that is, for $j\in\N\setminus \bigcup_{n=1}^{\infty}\operatorname{supp}(u_{n})$), set $\sigma_{j} = 1$. 
The $1$\nobreakdash-un\-con\-di\-tio\-na\-li\-ty of the unit vector basis $(e_j)_{j\in\N}$ for~$B_{p}$ means that we can define an isometric isomorphism $\Delta\in\mathscr{B}(B_p)$ by $\Delta x = \sum_{j=1}^{\infty}\sigma_{j}\langle x, e_j^*\rangle e_{j}$. Our choice of the sequence $(\sigma_j)_{j\in\N}$ implies that 
\begin{gather}\label{operatorontolpbasis:eq2} 
    \Delta(u_{n}) = |u_{n}| \qquad (n \in \N), 
\end{gather}
where we have used the standard notion of \emph{modulus} for an element of a Banach space with a $1$-un\-con\-di\-tional basis (justified by the fact that such a Banach space is a Banach lattice), that is, 
$\bigl\lvert\sum_{j=1}^{\infty}\alpha_{j}e_{j}\bigr\rvert = \sum_{j=1}^{\infty}\lvert\alpha_{j}\rvert e_{j}$.

For each $n \in \N$, take a Schreier chain~$\mathcal{C}_{n}$ contained in $\operatorname{supp}(u_{n})$ with $\beta_{p}(u_{n},\mathcal{C}_{n})=1$, and set $\mathcal{C} = \bigcup_{n \in \N}\mathcal{C}_{n}$. Defining $m_{0} = 0$ and $m_{n} = \sum_{k=1}^{n}\lvert\mathcal{C}_{k}\rvert$ for $n\in \N$,  we can enumerate~$\mathcal{C}_{n}$ as 
\begin{gather*}
    \mathcal{C}_{n} = \{F_{m_{n-1}+1}< F_{m_{n-1}+2} < \cdots < F_{m_{n}}\}.
\end{gather*}
Since $(u_n)_{n\in\N}$ is a block basic sequence, we have $F_{m_n}<F_{m_n+1}$ for $n\in\N$. Consequently \mbox{$\mathcal{C} = \{ F_1<F_2<\cdots\}$} is an infinite chain of successive Schreier sets, so it induces an operator $\Sigma_{\mathcal{C}}\in\mathscr{B}(B_{p},\ell_{p})$ of norm~$1$ by \Cref{Sigmaoperator}.

Set $D_{n} = \operatorname{span}(d_{j} : m_{n-1}< j \leq m_{n}) \subset \ell_{p}$, and let  $Q_{n}\in\mathscr{B}(\ell_{p},D_{n})$ be the basis projection on\-to~$D_n$; that is, $Q_nd_j = d_j$ for $m_{n-1}< j \leq m_{n}$ and $Q_nd_j = 0$ otherwise. Then we can define an isometric isomorphism $\Theta\in\mathscr{B}\bigl(\ell_{p},\bigl(\bigoplus_{n \in \N} D_{n}\bigr)_{\ell_{p}}\bigr)$ by $\Theta x = (Q_{n}x)_{n\in\N}$.

In view of~\eqref{operatorontolpbasis:eq2}, the vector $y_{n} = \Sigma_{\mathcal{C}} \Delta(u_{n})\in\ell_p$ satisfies
\begin{gather*}
    y_{n} = \Sigma_{\mathcal{C}}\lvert u_{n}\rvert = \sum_{j=m_{n-1}+1}^{m_{n}}\Bigl(\sum_{k\in F_{j}}\lvert\langle u_{n},e_{k}^{*}\rangle\rvert\Bigr)d_{j} \in D_{n}.
\end{gather*}
In particular, we have
\begin{gather*}
    \lv y_{n} \rv_{\ell_p}^{p}  = \sum_{j=m_{n-1}+1}^{m_{n}}\Bigl(\sum_{k\in F_{j}}\lvert\langle u_{n},e_{k}^{*}\rangle\rvert\Bigr)^p = \beta_{p}(u_{n},\mathcal{C}_{n})^{p} = 1,
\end{gather*}
so by the Hahn--Banach Theorem, we can take $f_{n}\in D_{n}^{*}$ such that $\langle y_{n},f_{n}\rangle = 1 = \lv f_{n}\rv_{\ell_p^*}$. This enables us to define an operator $\Gamma\in\mathscr{B}\bigl(\bigl(\bigoplus_{n \in\N}D_{n}\bigr)_{\ell_{p}},\ell_{p}\bigr)$ of norm~$1$ by $\Gamma(x_{n})_{n\in\N} = \bigl(\langle x_n,f_{n}\rangle\bigr)_{n\in \N}$.

Finally, we compose these operators to obtain an operator $V = \Gamma\Theta\Sigma_{\mathcal{C}}\Delta\in\mathscr{B}(B_p,\ell_p)$; that is,
\[ \begin{tikzcd}
	{B_{p}} && {B_{p}} && {\ell_{p}} && {\displaystyle{\Bigl(\bigoplus_{n\in\N}D_{n}\Bigr)_{\ell_{p}}}} && {\ell_{p}.}
	\arrow["\displaystyle{\Delta}", from=1-1, to=1-3]
	\arrow["\displaystyle{{\Sigma_{\mathcal{C}}}}", from=1-3, to=1-5]
	\arrow["\displaystyle{\Theta}", from=1-5, to=1-7] 
	\arrow["\displaystyle{\Gamma}", from=1-7, to=1-9]
\end{tikzcd} \]
Recalling that $y_{n} = \Sigma_{\mathcal{C}} \Delta(u_{n})\in D_n$ and then using the definitions of the operators~$\Theta$ and~$\Gamma$, we conclude that
\begin{gather*}
    Vu_n = \Gamma\Theta y_{n} = \bigl(\langle Q_jy_n,f_j\rangle\bigr)_{j\in\N} =    
    \langle y_{n},f_{n}\rangle d_{n} = d_{n}\qquad (n\in\N).  
\end{gather*}
In particular, since $u_n$ and $d_n$ are unit vectors, we have $1\le \lVert V\rVert\le \lVert\Gamma\rVert\,\lVert\Theta\rVert\,\lVert\Sigma_{\mathcal{C}}\rVert\,\lVert\Delta\rVert =1$, so~$V$ has norm~$1$.
\end{proof}

\begin{proof}[Proof of Theorem~{\normalfont{\ref{lpsaturation}}}]
By \Cref{C2afterC1} (applied with $X = E$ and $Y = D$), it suffices to show that for every constant $C>1$ and every block subspace $W= \overline{\operatorname{span}}\,(w_{n}:n \in \N)$ of~$E$, where  $(w_{n})_{n\in\N}$ is a block basic sequence of the unit vector basis 
for~$E$, there are operators $U\in\mathscr{B}(D,W)$ and $V\in\mathscr{B}(E,D)$ such that $V|_WU = I_D$ and $\lv U\rv\,\lv V\rv\le C$. 

\Cref{infinitynull} implies that we can find a normalized block basic sequence $(u_{n})_{n\in\N}$ of $(w_{n})_{n\in\N}$ such that $\lv u_{n} \rv_{\infty} \rightarrow 0$ as $n \rightarrow \infty$. By \Cref{dominatedbylp}, $(u_n)_{n\in\N}$ admits a subsequence $(u_{n_j})_{j\in\N}$ which is $C$-dominated by the unit vector basis $(d_j)_{j\in\N}$ for~$D$. Therefore we can define an operator $U\in\mathscr{B}(D,W)$ by $Ud_{j} = u_{n_j}$ for every $j\in \N$, and $\lv U \rv \leq C$. \Cref{operatorontolpbasis} shows that there exists an operator $V\in\mathscr{B}(E,D)$ of norm~$1$ such that $Vu_{n_j} = d_{j}$ for every $j\in \N$. It follows that $VUd_{j} = d_{j}$ for every $j \in \N$, so $V|_{W}U = I_{D}$, and $\lv U \rv \cdot \lv V \rv \leq C \cdot 1=C$, as required. 
\end{proof}

\begin{definition}
A Banach space~$X$ is \emph{sub\-projective} if every closed, in\-finite-di\-men\-sional subspace of~$X$ contains a~closed, infinite-dimensional subspace which is complemented in~$X$.    
\end{definition} 

\Cref{lpsaturation} implies that the Baernstein and Schreier spaces have this property.
We record this observation formally for later reference.

\begin{corollary}\label{C:subproj}
  The Baernstein spaces $B_{p}$ for $1<p<\infty$ and the Schreier spaces $S_p$ for $1\le p<\infty$ are subprojective.
\end{corollary}

\begin{remark}\label{R:Bpreflexive}
    Originally, Baernstein~\cite{AB} proved that the Banach space~$B_{2}$ is reflexive by verifying that the unit vector basis is a shrinking and boundedly complete basis for it and then appealing to a well-known theorem of James. We can now give an alternative proof of this result using \Cref{lpsaturation}, valid for any $1<p<\infty$: the fact that~$B_p$ is $\ell_p$-saturated implies that it does not contain any subspace isomorphic to either~$c_{0}$ or~$\ell_{1}$. Since~$B_p$ has an unconditional basis, it follows from another well-known theorem of James that~$B_{p}$ is reflexive.
\end{remark} 

With a small amount of extra effort, we can characterize the normalized block basic sequences of the unit vector basis for the Baernstein and Schreier spaces that admit a subsequence which is equivalent to the unit vector basis for~$\ell_p$ or~$c_0$, respectively. 

\begin{proposition}\label{P:charDsubseq}
    Let $(E,D) = (B_{p},\ell_{p})$ for some $1 < p < \infty$ or $(E,D) = (S_{p},c_{0})$ for some $1 \leq p < \infty$. The following conditions are equivalent for a normalized block basic sequence $(u_{n})_{n\in\N}$  of the unit vector basis for~$E\colon$
    \begin{enumerate}[label={\normalfont{(\alph*)}}]    
    \item\label{P:charDsubseq:a} $\inf_{n\in\N}\lVert u_n\rVert_\infty = 0;$
    \item\label{P:charDsubseq:b} $(u_{n})_{n\in\N}$ admits a subsequence which is $C$-equivalent to the unit vector basis for~$D$, for every constant $C>1;$
    \item\label{P:charDsubseq:c} $(u_{n})_{n\in\N}$ admits a subsequence which is dominated by the unit vector basis for~$D$. 
    \end{enumerate}
\end{proposition}

\begin{proof}
    To see that~\ref{P:charDsubseq:a} implies~\ref{P:charDsubseq:b}, suppose that $\inf_{n\in\N}\lVert u_n\rVert_\infty = 0$, and take $C>1$. By \Cref{dominatedbylp}, $(u_{n})_{n\in\N}$ admits a subsequence $(u_{n_j})_{j\in\N}$ that is $C$-dominated by $(d_j)_{j\in\N}$. On the other hand, $(u_{n_j})_{j\in\N}$ $1$-dominates $(d_j)_{j\in\N}$ because \Cref{operatorontolpbasis} shows that there is an operator $V\in\mathscr{B}(E,D)$ with $\lVert V\rVert =1$ such that $Vu_{n_j} = d_j$ for every $j\in\N$. Hence $(u_{n_j})_{j\in\N}$ and $(d_j)_{j\in\N}$ are $C$-equivalent. 

    The implication \ref{P:charDsubseq:b}$\Rightarrow$\ref{P:charDsubseq:c} is trivial.

    We complete the proof by proving that~\ref{P:charDsubseq:c} implies~\ref{P:charDsubseq:a}, arguing contrapositively. Suppose that $\delta := \inf_{n\in\N}\lVert u_n\rVert_\infty>0$, and take a subsequence $(u_{n_j})_{j\in\N}$ of   $(u_{n})_{n\in\N}$. To verify that $(d_j)_{j\in\N}$ does not dominate $(u_{n_j})_{j\in\N}$, it suffices to show that for every $C\ge 1$, there exists $k\in\N$ such that 
    \begin{equation}\label{P:charDsubseq:eq1}
        \biggl\lVert\sum_{j=k}^{2k-1} u_{n_j}\biggr\rVert_E > C \biggl\lVert\sum_{j=k}^{2k-1} d_j\biggr\rVert_D = \begin{cases} C\ &\text{for}\ D=c_0,\\ C k^{\frac1p}\ &\text{for}\ D=\ell_p.           
        \end{cases}
    \end{equation}
    Choose $k\in\N$ such that $k>(C/\delta)^p$ if $E=S_p$ and $k>(C/\delta)^{\frac{p}{p-1}}$ if $E=B_p$, and set $x= \sum_{j=k}^{2k-1} u_{n_j}\in E$. By hypothesis, we can find $m_j\in\operatorname{supp} (u_{n_j})$ such that $\lvert\langle u_{n_j}, e_{m_j}^*\rangle\rvert\ge\delta$ for each $j\in\{k,\ldots,2k-1\}$. Then $F= \{m_j : k\le j<2k\}$ is a Schreier set because $\lvert F\rvert = k\le m_k=\min F$. Hence we have
    \[ \lVert x\rVert_{S_p}\ge \mu_p(x,F) = \biggl(\sum_{j=k}^{2k-1} \lvert\langle x, e_{m_j}^*\rangle\rvert^p\biggr)^{\frac1p} = \biggl(\sum_{j=k}^{2k-1} \lvert\langle u_{n_j}, e_{m_j}^*\rangle\rvert^p\biggr)^{\frac1p}\ge k^{\frac1p}\delta >C \]
    and
    \[ \lVert x\rVert_{B_p}\ge \beta_p(x,\{F\}) = \sum_{j=k}^{2k-1} \lvert\langle x, e_{m_j}^*\rangle\rvert = \sum_{j=k}^{2k-1} \lvert\langle u_{n_j}, e_{m_j}^*\rangle\rvert\ge k\delta >Ck^{\frac1p}, \]
    where the final inequalities follow from the choice of~$k$ in both cases. This establishes~\eqref{P:charDsubseq:eq1}. 
\end{proof}

\section{An application to operator ideals: the proof of \Cref{thmsmalllargeideals}\ref{itemcompactss}}\label{S:3}

\noindent The main purpose of this short section is to use \Cref{lpsaturation} to  identify the ideals of strictly singular, inessential and weak\-ly compact operators on the Baernstein and Schreier spaces. We begin by recalling the formal definitions of these ideals, as well as some other standard notions that we require. 

\begin{definition}
An operator $T\in\mathscr{B}(X,Y)$ between Banach spaces~$X$ and~$Y$ is:
\begin{itemize}
\item \emph{strictly singular} if the restriction of~$T$ to~$W$ is not an isomorphic embedding for any infinite-dimensional subspace $W$ of $X$,
\item \emph{inessential} if $I_{X} + UT$ is a Fredholm operator (meaning that its kernel is finite-di\-men\-sional and its range has finite codimension in~$X$) for every operator $U\in\mathscr{B}(Y,X)$,
\item \emph{weakly compact} if the image under~$T$ of the unit ball in~$X$ is relatively weakly compact,
\item \emph{unconditionally converging} if the series $\sum_{n=1}^{\infty}Tx_{n}$ converges unconditionally in norm for every series $\sum_{n=1}^{\infty}x_{n}$ in~$X$ which is weakly unconditionally Cauchy in the sense that the series $\sum_{n=1}^{\infty}\langle x_{n}, f\rangle$ converges absolutely for every functional $f\in X^{*}$.
\end{itemize}
Furthermore, we say that $T$ \emph{fixes a copy} of a Banach space~$Z$ if there is an operator $V\in\mathscr{B}(Z,X)$ such that the composition $TV$ is an isomorphic embedding. 
\end{definition}
\noindent We write $\mathscr{S}(X,Y)$, $\mathscr{E}(X,Y)$ and $\mathscr{W}(X,Y)$ for the sets of strictly singular, inessential and weakly compact operators between~$X$ and~$Y$, respectively, with the usual convention that   $\mathscr{S}(X) = \mathscr{S}(X,X)$, \textit{etc.} It is well known that~$\mathscr{S}$, $\mathscr{E}$ and~$\mathscr{W}$ are closed operator ideals in the sense of Pietsch.

The Banach--Alaoglu Theorem implies that every operator defined on a reflexive Banach space is weakly compact, so $\mathscr{W}(B_p) = \mathscr{B}(B_p)$ for every $1<p<\infty$. 

Pfaffenberger~\cite{Pfaf} has shown that $\mathscr{S}(X) = \mathscr{E}(X)$ for every subprojective Banach space~$X$. Hence, in view of \Cref{C:subproj}, we obtain: 

\begin{proposition}\label{P:SS=iness}
    Let $E= B_p$ for some $1<p<\infty$ or $E=S_p$ for some $1\le p<\infty$. Then \[ \mathscr{S}(E) = \mathscr{E}(E). \]
\end{proposition}

Consequently, to complete the proof of \Cref{thmsmalllargeideals}\ref{itemcompactss}, it remains only to show that $\mathscr{W}(S_p) = \mathscr{S}(S_p)$ for $1\le p<\infty$. This will follow from our next, considerably more general, result, which applies to~$X=S_p$ by \Cref{lpsaturation} and the fact that~$S_p$ has an unconditional basis. 

\begin{proposition}\label{strictlysingularweaklycompact}
Let $T\in\mathscr{B}(X,Y)$ be an operator, where~$X$ is a $c_0$-saturated Banach space that em\-beds into a Banach space with an unconditional basis and~$Y$ is a separable Banach space. The following conditions are equivalent:
\begin{enumerate}[label={\normalfont{(\alph*)}}]    
\item\label{strictlysingularweaklycompact:i} $T$ is strictly singular;
\item\label{strictlysingularweaklycompact:ii} $T$ is inessential;
\item\label{strictlysingularweaklycompact:vi} the identity operator on~$c_0$ does not factor through~$T$ in the sense that there are no ope\-ra\-tors $U\in\mathscr{B}(Y,c_0)$ and $V\in\mathscr{B}(c_0,X)$ such that $UTV=I_{c_0};$
\item\label{strictlysingularweaklycompact:v} $T$ does not fix a copy of~$c_0;$
\item\label{strictlysingularweaklycompact:iv} $T$ is unconditionally converging;
\item\label{strictlysingularweaklycompact:iii} $T$ is weakly compact.
\end{enumerate}
\end{proposition}

The proof relies on a classical result of Pe\l{}czy\'{n}ski \cite[Proposition~9, $1^\circ$]{P}.

\begin{theorem}\label{T:Pelc}
Let~$X$ be a Banach space which embeds into a Banach space with an unconditional basis, and suppose that~$X$ does not contain any subspace which is isomorphic to~$\ell_1$. Then every unconditionally converging operator from~$X$ into a Banach space is weakly compact.
\end{theorem}

\begin{proof}[Proof of Proposition~{\normalfont{\ref{strictlysingularweaklycompact}}}]
The proof has two parts. In part~(i), we show that conditions~\ref{strictlysingularweaklycompact:i}--\ref{strictlysingularweaklycompact:v} are equivalent, while part~(ii) contains the proof that conditions~\ref{strictlysingularweaklycompact:v}--\ref{strictlysingularweaklycompact:iii} are equivalent.

(i). The implication \ref{strictlysingularweaklycompact:i}$\Rightarrow$\ref{strictlysingularweaklycompact:ii} is always true, and~\ref{strictlysingularweaklycompact:ii} implies~\ref{strictlysingularweaklycompact:vi}  because the identity operator on an infinite-dimensional Banach space cannot be inessential. 

We prove that~\ref{strictlysingularweaklycompact:vi} implies~\ref{strictlysingularweaklycompact:v} by contraposition. Suppose that we can find an operator $V\in\mathscr{B}(c_0,X)$ such that~$TV$ is an isomorphic embedding. Then~$TV[c_0]$ is isomorphic to~$c_0$, so it is complemented in~$Y$ by Sobczyk's Theorem. Therefore~$TV$ has a left inverse $U\in\mathscr{B}(Y,c_0)$; that is, $UTV=I_{c_0}$, as desired. 

Finally, \ref{strictlysingularweaklycompact:v} implies~\ref{strictlysingularweaklycompact:i} because~$X$ is $c_0$-saturated. 

(ii). As observed in \cite[Exercise~8(i), page~54]{JD}, conditions~\ref{strictlysingularweaklycompact:v} and~\ref{strictlysingularweaklycompact:iv} are equivalent in general (that is, without any restrictions on the Banach spaces~$X$ and~$Y$). The hypothesis that~$X$ is $c_0$-saturated ensures that~$\ell_1$ does not embed into~$X$, so 
\Cref{T:Pelc} shows that~\ref{strictlysingularweaklycompact:iv} im\-plies~\ref{strictlysingularweaklycompact:iii}. Finally, the implication \ref{strictlysingularweaklycompact:iii}$\Rightarrow$\ref{strictlysingularweaklycompact:v} follows from the fact that~$c_0$ is not reflexive. 
\end{proof}

\begin{remark} The condition that the codomain~$Y$ of the operator~$T$ in \Cref{strictlysingularweaklycompact} is separable cannot be dropped because the embedding of~$c_0$ into~$\ell_\infty$ is an inessential opeator which obviously fixes a copy of~$c_0$.     
\end{remark}

\section{The Gasparis--Leung index and its applications}\label{S:4}
\noindent The aim of this section is to establish counterparts for the Baernstein and Schreier spaces of some important  technical results of Gasparis and Leung~\cite{GL}. They introduced a numerical index for each $n\in\N$ and every pair $M,N$ of infinite subsets of~$\N$ which characterizes when the subspaces spanned by the infinite subsequences of the unit vector basis for~$X[\mathcal{S}_n]$ corresponding to~$M$ and~$N$ are isomorphic. 
As we mentioned in the introduction, it turns out that this index, for $n=1$, also works for the Baernstein and Schreier spaces. 

In order to define it, we must first introduce the \emph{Schreier covering number} of a set $A \in [\N]^{< \infty}$, which according to \cite[Definition~3.1]{GL} and using our notation from \Cref{S:2} is
\begin{equation}\label{Eq:SCN}
\tau_1(A) = \begin{cases} 0\ &\text{if}\ A=\emptyset,\\
    \min\Bigl\{\lvert\mathcal{C}\rvert : \mathcal{C}\in\SC,\, A\subseteq\bigcup\mathcal{C} \Bigr\}\ &\text{otherwise.}    
\end{cases} 
\end{equation} 
Unpacking the somewhat condensed notation for $A\ne\emptyset$, we can restate this definition as  
\[ \tau_1(A) = \min\biggl\{ m\in \mathbb{N} : A \subseteq \bigcup_{j=1}^{m}F_{j},\ \text{where}\ F_1,\ldots,F_m\in\mathcal{S}_1\ \text{and}\ F_1<F_2<\cdots<F_m\biggr\}. \]
Furthermore, as observed in \cite[Remark 4.2]{BKL}, we can refine it as follows. Let $m\in\N$. Then $\tau_1(A) = m$ if and only if there is a Schreier chain \mbox{$\{F_{1} < \cdots < F_{m}\}$} such that $A = \bigcup_{j=1}^{m}F_{j}$ and $F_{1},\ldots,F_{m-1}$ are maximal Schreier sets; it is important to note that~$F_{m}$ need not be maximal. 

As in \cite{BKL}, for a set $M = \{m_1<m_2<\cdots\}\in [\mathbb{N}]$ and $J \subseteq \mathbb{N}$, we define
\[ M(J) = \{m_{j} : j \in J\}. \]
This piece of notation enables us to state \cite[Definition~3.3]{GL} in the  following compact form. For  $M,N\in[\mathbb{N}]$, the \emph{Gas\-pa\-ris--Leung index} of~$M$ with respect to~$N$ is 
\begin{equation}\label{Eq:GLindex}
\gl(M,N) = \sup\{\tau_1(M(J)) : J \in [\mathbb{N}]^{< \infty},\, N(J) \in \mathcal{S}_{1}\}.
\end{equation}
Gasparis and Leung denoted this quantity~$d_1(M,N)$. We have chosen the more distinctive sym\-bol~$\gl(M,N)$ in their honour, noting that the Greek spelling of ``Gasparis'' begins with the letter~$\Gamma$.

Before we state the first main result of this section, let us introduce a piece of notation that we shall use frequently. Given a Banach space~$X$ with a basis $(x_n)_{n\in\N}$, we set
\begin{equation}\label{eq:X_N}
X_N = \overline{\operatorname{span}}\, (x_n : n\in N)\qquad (N\subseteq\N).    
\end{equation} 
Further, recall from \Cref{S:2} that we write $\lVert x\rVert_\infty = \sup_{n\in\N}\lvert \langle x,e_{n}^*\rangle\rvert$ for $x\in B_p$ or $x\in S_p$ in line with standard usage.

\begin{theorem}\label{GLThm}
    Let $E= B_p$ for some $1<p<\infty$ or $E=S_p$ for some $1\le p<\infty$, equipped with the unit vector basis $(e_{n})_{n\in\N}$. 
The following conditions are equivalent for $M,N\in[\N]\colon$
        \begin{enumerate}[label={\normalfont{(\alph*)}}]
        \item\label{GLThm1ii} The Gasparis--Leung index $\gl(M,N)$ is finite.
        \item\label{GLThm1i} The basic sequence $(e_{m})_{m\in M}$ dominates $(e_{n})_{n\in N}$.
        \item\label{GLThm1iii} There exists an operator $T\in\mathscr{B}(E_M,E_N)$ for which $\inf_{m\in M}\lVert Te_{m}\rVert_\infty > 0$.
    \end{enumerate}
\end{theorem}

We begin with a quantitative version of the implication \ref{GLThm1ii}$\Rightarrow$\ref{GLThm1i}.
\begin{lemma}\label{L:GLthmAimpliesB}
    Let $E = B_p$ for some $1<p<\infty$ or $E = S_p$ for some $1\le p<\infty$, take $M,N\in[\N]$ for which $\gl(M,N)<\infty$, and define
    \begin{equation}\label{L:GLthmAimpliesB:eq0}  
    C = \begin{cases} 
        \gl(M,N)\ &\text{for}\ E = B_p,\\ 
        \gl(M,N)^{\frac1p}\ &\text{for}\ E = S_p.
    \end{cases} 
    \end{equation}
  Then the basic sequence $(e_m)_{m\in M}$ $C$-dominates $(e_n)_{n\in N}$.
\end{lemma}

\begin{proof} Enumerate~$M$ and~$N$ as $M = \{m_1<m_2<\cdots\}$ and $N=\{n_1<n_2<\cdots\}$, respectively. 
Our aim is to show that $\lVert y\rVert_{E}\le C\lVert x\rVert_{E}$ whenever $x = \sum_{j=1}^k\alpha_j e_{m_j}$ and $y = \sum_{j=1}^k\alpha_j e_{n_j}$ for some $k\in\N$ and some $\alpha_1,\ldots,\alpha_k\in\mathbb{K}$. We may of course suppose that 
$\alpha_1,\ldots,\alpha_k$ are not all~$0$. 

We consider the Baernstein and Schreier spaces separately, but emphasize that the proofs follow similar strategies, originating in the proof of \cite[Lemma~3.4]{GL}. 
For readability, we begin with the easier case, which is $E = S_p$.  
 Since~$y$ is finitely supported, we can choose a Schreier set~$F$ such that $\lVert y\rVert_{S_p} = \mu_p(y,F)$ and $F\subseteq \operatorname{supp} y\subseteq \{n_i : 1\le i\le k\}$. Take $J\subseteq\{1,\ldots, k\}$ such that $N(J) = F$. Then $\tau_1(M(J))\leq\gl(M,N) = C^p$ by~\eqref{Eq:GLindex} and~\eqref{L:GLthmAimpliesB:eq0}, so there is a Schreier chain $\{G_{1} < \cdots < G_{C^p}\}$ such that $M(J)\subseteq\bigcup_{i=1}^{C^p}G_{i}$ by~\eqref{Eq:SCN}. Set $K_{i} = \{ j\in J : m_j\in G_{i}\}$ for $i\in\{1,\ldots,C^p\}$. Then we have 
$J = \bigcup_{i=1}^{C^p}K_{i}$ and $K_h\cap K_i=\emptyset$ for $h\ne i$, from which we deduce that
\begin{gather*}
 \lVert y\rVert_{S_p}^p = \mu_p(y,F)^p = \sum_{j\in J}\lvert\alpha_{j}\rvert^{p} = \sum_{i=1}^{C^p}\sum_{j\in K_{i}}\lvert \alpha_{j}\rvert^{p}  = \sum_{i=1}^{C^p}\mu_{p}(x, G_{i})^{p}\leq C^p\lVert x\rVert_{S_{p}}^{p}, 
\end{gather*}
where the final inequality follows from the fact that $G_1,\ldots,G_{C^p}\in\mathcal{S}_1$. 

Having completed the proof for $E = S_p$, we turn our attention to $E = B_p$. 
We begin in the same way as above: using that~$y$ is finitely supported, we can find a Schreier chain $\mathcal{C} = \{F_{1} < \cdots < F_{t}\}$ such that $\lVert y\rVert_{B_p} = \beta_p(y,\mathcal{C})$ and $\bigcup_{r=1}^t F_r\subseteq\operatorname{supp} y\subseteq\{n_i : 1\le i\le k\}$. 

Fix $r\in\{1,\ldots,t\}$, and choose $J_r\subseteq\{1,\ldots,k\}$ such that $N(J_{r}) = F_{r}\in\mathcal{S}_1$.  Then we have $\tau_1(M(J_{r}))\leq \gl(M,N) = C$ by~\eqref{Eq:GLindex} and~\eqref{L:GLthmAimpliesB:eq0}, so we can find a Schreier chain \mbox{$\{G_{1}^{r} <\cdots < G_{C}^{r}\}$} such that $M(J_{r})\subseteq \bigcup_{i=1}^{C}G_{i}^{r}$ by~\eqref{Eq:SCN}.
Set 
\[ K_{i}^{r} = \{j\in J_{r} : m_{j} \in G_{i}^{r}\}\qquad \text{and}\qquad \gamma_i^r = \sum_{j\in K_i^r}\lvert\alpha_j\rvert\qquad (i\in\{1,\ldots,C\}),  \]
and choose $\iota(r)\in\{1,\ldots,C\}$ for which $\gamma_{\iota(r)}^r = \max\{\gamma_i^r : 1\le i\le C\}$. Since $J_{r} = \bigcup_{i=1}^{C} K_{i}^{r}$ and $K_{h}^{r} \cap K_{i}^{r} = \emptyset$ whenever $h\ne i$, we have 
\begin{equation*}
\sum_{j\in J_{r}}\lvert\alpha_j\rvert = \sum_{i=1}^{C} \gamma_i^r\le C\gamma_{\iota(r)}^r = C\sum_{j\in K_{\iota(r)}^r}\lvert\alpha_j\rvert.
\end{equation*}
Furthermore, $\mathcal{D} = \{ M(K_{\iota(r)}^{r}) : 1\le r\le t\}$ is a Schreier chain, as the following two facts show:
\begin{itemize}
\item \mbox{$M(K_{\iota(r)}^{r}) \in \mathcal{S}_{1}$} for each $1\le r\le t$ because $M(K_{\iota(r)}^{r}) \subseteq G_{\iota(r)}^{r} \in \mathcal{S}_{1}$, and
\item the sets $M(K_{\iota(1)}^{1}),M(K_{\iota(2)}^{2}),\ldots,M(K_{\iota(t)}^{t})$ are successive because $M(K_{\iota(r)}^r)\subseteq M(J_r)$ for each $1\le r\le t$ and $M(J_1)<M(J_2)<\cdots<M(J_t)$.  
\end{itemize}
In conclusion, we have 
\[ \lVert y\rVert_{B_p}^p = \beta_p(y,\mathcal{C})^p = \sum_{r=1}^t\biggl(\sum_{j\in J_r}\lvert\alpha_j\rvert\biggr)^p\le C^p\sum_{r=1}^t\biggl(\sum_{j\in K_{\iota(r)}^r}\lvert\alpha_j\rvert\biggr)^p = C^p\beta_p(x,\mathcal{D})^p\le C^p\lVert x\rVert_{B_p}^p. \qedhere  \]
\end{proof}

\begin{proof}[Proof of Theorem~{\normalfont{\ref{GLThm},}} \ref{GLThm1ii}$\Rightarrow$\ref{GLThm1i}$\Rightarrow$\ref{GLThm1iii}] \Cref{L:GLthmAimpliesB} shows that~\ref{GLThm1ii} implies~\ref{GLThm1i}. 

To see that~\ref{GLThm1i} implies~\ref{GLThm1iii}, suppose that $(e_{m})_{m\in M}$ dominates $(e_{n})_{n\in N}$. By definition, this means that the linear map $T\colon \operatorname{span}(e_{m}:m\in M) \rightarrow E_N$ determined by $Te_{m_{j}} = e_{n_{j}}$ for every $j \in \N$ is bounded, where $M = \{m_1<m_2<\cdots\}$ and \mbox{$N=\{n_1<n_2<\cdots\}$} are the increasing enumera\-tions. Therefore~$T$ extends uniquely to an operator in~$\mathscr{B}(E_M,E_N)$, also denoted~$T$, which satisfies $\langle Te_{m_{j}},e_{n_{k}}^{*}\rangle = \delta_{j,k}$ for every $j,k\in \N$. Hence $\inf_{m\in M}\lVert Te_{m}\rVert_\infty = 1 > 0$.
\end{proof}

It remains to prove the implication \ref{GLThm1iii}$\Rightarrow$\ref{GLThm1ii} in \Cref{GLThm}. For this, we follow the approach of~\cite{GL} closely, although we can shorten certain steps because we consider only the first Schreier family~$\mathcal{S}_1$. 

We begin by generalizing \cite[Proposition~3.13]{GL}, which Gasparis and Leung established for the higher-order Schreier spaces $X[\mathcal{S}_\xi]$ for $\xi<\omega_1$. However, as we shall show, it applies to a much larger class of Banach spaces, including the Baernstein and Schreier spaces that we are investigating. We provide a detailed proof for the reader's convenience.  

\begin{lemma}\label{1.1step1BpV2}
Let~$X$ and~$Y$ be Banach spaces with unconditional, normalized bases $(x_{n})_{n\in\N}$ and $(y_n)_{n\in\N}$, respectively, and suppose that the basis~$(x_n)_{n\in\N}$ for~$X$ is shrinking.
The following con\-di\-tions are equivalent: 
\begin{enumerate}[label={\normalfont{(\alph*)}}]
    \item\label{1.1step1BpV2:a} There is an operator $T\in\mathscr{B}(X,Y)$ for which 
\begin{equation}\label{1.1step1BpV2:eq1}
  \inf_{k\in\N}\sup_{j\in\N}\lvert\langle Tx_{k},y_{j}^{*}\rangle\rvert > 0.
\end{equation}
\item\label{1.1step1BpV2:b} There is an operator $U\in\mathscr{B}(X,Y)$ for which $Ux_k\in \{ y_j : j\in\N\}$ for every $k\in\N$. 
\end{enumerate}
\end{lemma}

It is easy to see that~\ref{1.1step1BpV2:b} implies~\ref{1.1step1BpV2:a}. The proof of the converse relies on a careful analysis of the matrix $(T_{j,k})_{j,k\in\N}$ associated with the operator~$T$. We recall the standard definition of this matrix: for $j,k\in \N$, the $(j,k)^{\text{th}}$ coefficient of the matrix associated with an operator \mbox{$T\in\mathscr{B}(X,Y)$} between Banach spaces~$X$ and~$Y$ with bases~$(x_{n})_{n\in\N}$ and $(y_{n})_{n\in\N}$, respectively, is 
\begin{equation}\label{Defn:MatrixofOperator}
    T_{j,k} = \langle Tx_k,y_j^*\rangle;
\end{equation}
that is, the $k^{\text{th}}$ column of the matrix associated with~$T$ contains the coordinates with respect to the basis~$(y_j)_{j\in\N}$ of the image under~$T$ of the $k^{\text{th}}$ basis vector~$x_k$. Dualizing, we have $T_{j,k} = \langle x_k,T^*y_j^*\rangle$, so if the basis $(x_n)_{n\in\N}$ for~$X$ is shrinking, then the $j^{\text{th}}$ row of the matrix associated with~$T$ contains the coordinates with respect to the basis $(x^*_k)_{k\in\N}$ for~$X^*$ of the image under~$T^*$ of the $j^{\text{th}}$ coordinate functional~$y_j^*$.

In the proof of \Cref{1.1step1BpV2}, we require a variant of a result due to Tong~\cite{Tong}. It involves the following notion: a  matrix $\Gamma = (\gamma_{j,k})_{j,k\in\N}$ is a \emph{block diagonal} of a matrix $A = (\alpha_{j,k})_{j,k\in\N}$ if there are increasing sequences $0\le r_{1} < r_{2} < \cdots$ and $0\le s_{1} < s_{2} < \cdots$ of integers for which
\begin{equation*}
\gamma_{j,k} = \begin{cases}
\alpha_{j,k}\ &\text{if}\ (j,k)\in \bigcup_{i=1}^{\infty}(r_{i},r_{i+1}] \times (s_{i},s_{i+1}]\\
0 &  \text{otherwise} 
\end{cases}\qquad (j,k\in\N).    
\end{equation*}

\begin{lemma}\label{lindenlemma}
Let $T\in\mathscr{B}(X,Y)$ be an operator between Banach spaces~$X$ and~$Y$ with un\-con\-di\-tional bases~$(x_{n})_{n\in\N}$ and $(y_{n})_{n\in\N}$, respectively, and suppose that $\Gamma = (\gamma_{j,k})_{j,k\in\N}$ is a block di\-ag\-o\-nal of the matrix associated with~$T$. Then there is an operator $R\in\mathscr{B}(X,Y)$ whose matrix is~$\Gamma;$ that is, 
\[ \langle Rx_k, y_j^*\rangle = \gamma_{j,k}\qquad (j,k\in\N). \] 
\end{lemma}

\begin{proof} As already mentioned, Tong proved a similar result in~\cite{Tong}, using very different terminology. A simple proof of the above statement is outlined in the first remark after \cite[Proposition~1.c.8]{L}. Note, however, that the definition stated in the text above  \cite[Proposition~1.c.8]{L} of the matrix associated with an operator~$T\in\mathscr{B}(X,Y)$ produces the transpose of the matrix given by~\eqref{Defn:MatrixofOperator}. Fortunately this difference does not matter, as the transpose of a block diagonal is again a block diagonal of the transposed matrix.
\end{proof}

\begin{proof}[Proof of Lemma~{\normalfont{\ref{1.1step1BpV2}}}] To see that~\ref{1.1step1BpV2:b} implies~\ref{1.1step1BpV2:a}, suppose that $U\in\mathscr{B}(X,Y)$ is an operator for which $Ux_k\in \{ y_j : j\in\N\}$ for every $k\in\N$. Then $\sup_{j\in\N}\lvert\langle Ux_{k},y_{j}^{*}\rangle\rvert = 1$ for every  $k\in\N$, so $T=U$ satisfies~\eqref{1.1step1BpV2:eq1}. 

We prove that~\ref{1.1step1BpV2:a} implies~\ref{1.1step1BpV2:b} by expanding on the approach Gasparis and Leung took in their proof of \cite[Proposition~3.13]{GL}. 
In view of~\eqref{1.1step1BpV2:eq1} and~\eqref{Defn:MatrixofOperator}, we can choose $\delta>0$ such that, for every $k\in\N$, $\lvert T_{j,k}\rvert\ge\delta$ for some $j\in\N$. This allows us to define a map $\psi\colon\N\to\N$ by \[ \psi(k) = \min\{j\in\N : \lvert T_{j,k}\rvert\ge\delta\}. \] Take $j\in\N$. Since the basis~$(x_n)_{n\in\N}$ for~$X$ is shrinking, the series $\sum_{k=1}^\infty T_{j,k}x_k^*$ is convergent with sum~$T^*y_j^*$, as explained in the text below~\eqref{Defn:MatrixofOperator}. It follows that $\lvert T_{j,k}\rvert\to0$ as $k\to\infty$, so every natural number has finite (possibly empty) pre-image under~$\psi$. There\-fore~$\psi$ has infinite image; let $\psi(\N) = \{n_1<n_2<\cdots\}$ be its increasing enumeration. Then $\{ \psi^{-1}(n_j) : j\in\N\}$ partitions~$\N$ into non-empty, finite, disjoint sets, so there is a unique permutation $\rho\colon\N\to\N$ such that 
\[ \rho(k)<\rho(m)\ \Longleftrightarrow\ \begin{cases} \psi(k)<\psi(m),\ \text{or}\\ \psi(k)=\psi(m)\ \text{and}\ k<m\end{cases}\qquad (k,m\in\N). \] 
In more concrete terms, we can define~$\rho$  as follows. Set $s_0 = 0$ and $s_j = \sum_{i=1}^j\lvert\psi^{-1}(n_i)\rvert$ for $j\in\N$. Then each $k\in\N$ belongs to the interval $(s_{j-1},s_j]$ for a unique~$j\in\N$, and $\rho(k)$ is the $(k-s_{j-1})^{\text{th}}$ smallest element of the set~$\psi^{-1}(n_j)$. In particular, we have $\psi(\rho(k)) = n_j$, so the definition of~$\psi$ implies that $\lvert T_{n_j,\rho(k)}\rvert\ge\delta$, and therefore we can define a diagonal operator $\Delta\in\mathscr{B}(X)$ of norm at most~$K/\delta$ by $\Delta x_{m} = T_{n_j,m}^{-1}x_{m}$ for each $m\in\N$, where~$j\in\N$ is chosen such that $s_{j-1}<\rho^{-1}(m)\le s_j$
and~$K$ denotes the unconditional basis constant of~$(x_n)_{n\in\N}$. 

The unconditionality of the basis~$(x_n)_{n\in\N}$ means that any reordering of it is also a basis for~$X$.
Hence, viewing the composite operator~$P_{\psi(\N)}T\Delta$ as a map from~$X$ to $Y_{\psi(\N)} = \overline{\operatorname{span}}\, (y_{n_j} : j\in\N)$, we may consider its matrix with respect to the bases~$(x_{\rho(k)})_{k\in\N}$ for~$X$  and~$(y_{n_j})_{j\in\N}$ for~$Y_{\psi(\N)}$. Suppose that $j,k\in\N$ satisfy $s_{j-1}<k\le s_j$. Then we have 
\[ (P_{\psi(\N)}T\Delta)_{j,k} = \langle  P_{\psi(\N)}T\Delta x_{\rho(k)},y_{n_j}^*\rangle = \frac{\langle  Tx_{\rho(k)},P_{\psi(\N)}^*y_{n_j}^*\rangle}{T_{n_j,\rho(k)}} = \frac{\langle Tx_{\rho(k)},y_{n_j}^*\rangle}{T_{n_j,\rho(k)}} = 1,  \] 
so the matrix $\Gamma = (\gamma_{j,k})_{j,k\in\N}$ defined by
\begin{equation}\label{1.1step1BpV2:eq2}
\gamma_{j,k} = \begin{cases} 1\ &\text{if}\ s_{j-1}<k\le s_j\\ 0\ &\text{otherwise} \end{cases}\qquad (j,k\in\N) 
\end{equation}
is  a block diagonal of the matrix associated with the operator~$P_{\psi(\N)}T\Delta$.  
\Cref{lindenlemma} implies that there is an operator $R\in\mathscr{B}(X,Y_{\psi(\N)})$ whose matrix is~$\Gamma$; that is, $\langle Rx_{\rho(k)},y_{n_j}^*\rangle = \gamma_{j,k}$ for every $j,k\in\N$. In view of~\eqref{1.1step1BpV2:eq2}, this means that $Rx_{\rho(k)} = y_{n_j}$ for every $k\in\N$, where $j\in\N$ is the unique number such that $s_{j-1}<k\le s_j$. Hence, writing $J\colon Y_{\psi(\N)}\to Y$ for the inclusion map, we obtain an operator $U=JR\in\mathscr{B}(X,Y)$ which satisfies $Ux_m = Rx_m\in \{y_n : n\in\N\}$ for every $m\in\N$, as required.
\end{proof}

\begin{lemma}\label{L:1Dec2023}
    Let $x = \sum_{j=1}^k \alpha_j e_{m_j}$, where $k\in\N$, $\alpha_1,\ldots,\alpha_k\in[0,\infty)$, and $m_1,\ldots,m_k\in\N$ are (not necessarily distinct) numbers for which $\{ m_1,\ldots,m_k\}\in\mathcal{S}_1$. Then $\lVert x\rVert_E = \sum_{j=1}^k\alpha_j$ for $E=S_1$ and $E=B_p$, while $\lVert x\rVert_{S_p}\ge\bigl(\sum_{j=1}^k\alpha_j^p\bigr)^{1/p}$ for $1<p<\infty$.
\end{lemma}

\begin{proof}
   Take $J\subseteq\{1,\ldots,k\}$ such that $\{m_j : j\in J\} = \{ m_1,\ldots,m_k\}$ and $m_i\ne m_j$ for  distinct $i,j\in J$, and set $K_j = \{ i\in\{1,\ldots,k\} : m_i=m_j\}$ for each $j\in J$. Then $\{ K_j : j\in\ J\}$ partitions $\{1,\ldots,k\}$, and we have $x= \sum_{j\in J}\bigl(\sum_{i\in K_j}\alpha_i\bigr) e_{m_j}$.
   Since $\{m_j : j\in J\}$ is  a Schreier set, we conclude that $\lVert x\rVert_E = \sum_{j\in J}\bigl(\sum_{i\in K_j}\alpha_i\bigr)= \sum_{j=1}^k\alpha_j$ for $E=S_1$ and $E=B_p$,    while 
   \[ \lVert x\rVert_{S_p}^p = \sum_{j\in J}\Bigl(\sum_{i\in K_j}\alpha_i\Bigr)^p\ge \sum_{j\in J} \sum_{i\in K_j}\alpha_i^p = \sum_{j=1}^k\alpha_j^p  \]
for $1<p<\infty$,  where the inequality follows from the fact that the $\ell_1$-norm dominates the $\ell_p$-norm.
\end{proof}

\begin{lemma}\label{L:1.1step2Bp} Let $E= B_p$ for some $1<p<\infty$ or $E=S_p$ for some $1\le p<\infty$, let $M\in [\N]$, and suppose that $\theta\colon M\to\N$ is a map for which the linear map $U\colon \operatorname{span} (e_m : m\in M)\to E$ determined by $U e_m = e_{\theta(m)}$ for $m\in M$ is bounded. Then
\begin{equation}\label{L:1.1step2Bp:eq1}
\sup\{\lvert\theta^{-1}(n)\rvert : n\in\N\}<\infty\qquad\text{and}\qquad  \sup \{\tau_1(\theta^{-1}(F)) : F\in\mathcal{S}_1\}<\infty.   
\end{equation}
\end{lemma}

\begin{proof} The hypothesis means that $U$ extends uniquely to an operator in~$\mathscr{B}(E_M,E)$, also denoted~$U$. 
We begin by showing that the second supremum in~\eqref{L:1.1step2Bp:eq1} is finite. Note that this will in\-clude showing that the pre-image under~$\theta$ of every Schreier set~$F$ is finite, as otherwise $\tau_1(\theta^{-1}(F))$ is not defined. 
Our strategy is as follows: given $F\in\mathcal{S}_1$, we take $m\in\N$ for which $\theta^{-1}(F)$ contains a chain $\{G_1<G_2<\cdots<G_m\}$ of maximal Schreier sets. As we shall verify below, $m$~is then dominated by a constant times a power of the norm of the operator~$U$; that is, $m\le C\lVert U\rVert^t$ for some constants $C,t\in(0,\infty)$ that will depend only on~$p$. This will give the desired conclusion because (i)~if~$\theta^{-1}(F)$ were infinite, it would contain arbitrarily long chains of maximal Schreier sets, contradicting the uniform bound on~$m$; (ii)~we can therefore use the characterization of~$\tau_1$ stated in the paragraph below its definition~\eqref{Eq:SCN} to deduce that $\tau_1(\theta^{-1}(F))\le C\lVert U\rVert^t+1$. Since the right-hand side of this inequality is independent of~$F$, it provides an upper bound on the second supremum in~\eqref{L:1.1step2Bp:eq1}.

It remains to establish the inequality $m\le C\lVert U\rVert^t$. We consider the two types of spaces separately. For $E=B_p$, set $x=\sum_{k=1}^m \lvert G_k\rvert^{-1}\sum_{j\in G_k} e_j$, 
and recall from~\eqref{estimate2baernstein} that \mbox{$\lVert x\rVert_{B_p}^p\le 2^p m$}. Consequently, we have 
\begin{equation}\label{L:1.1step2Bp:eq3} 2^p m\lVert U\rVert^p\ge \lVert Ux\rVert_{B_p}^p = \biggl\lVert \sum_{k=1}^m \frac{1}{\lvert G_k\rvert}\sum_{j\in G_k} e_{\theta(j)}\biggr\rVert_{B_p}^p = \biggl(\sum_{k=1}^m \frac{\lvert G_k\rvert}{\lvert G_k\rvert}\biggr)^p = m^p, \end{equation}
where the penultimate step follows from \Cref{L:1Dec2023} and the fact that $\bigcup_{k=1}^m\theta(G_k)$ is a Schreier set because it is contained in $F\in\mathcal{S}_1$. Rearranging~\eqref{L:1.1step2Bp:eq3}, we obtain $m\le (2\lVert U\rVert)^{\frac{p}{p-1}}$, which provides an upper bound on~$m$ of the desired form for $C=2^{\frac{p}{p-1}}$ and $t=p/(p-1)$.

The argument for $E=S_p$ is very similar, except that we use the vector \[ x=\sum_{k=1}^m \frac{1}{\lvert G_k\rvert^{\frac1p}}\sum_{j\in G_k} e_j. \] It has $S_p$-norm at most~$2^{\frac1p}$ by~\eqref{estimate2schreier}, so  
\[ 2\lVert U\rVert^p\ge \lVert U x\rVert_{S_p}^p = \biggl\lVert \sum_{k=1}^m \frac{1}{\lvert G_k\rvert^{\frac1p}}\sum_{j\in G_k} e_{\theta(j)}\biggr\rVert_{S_p}^p\ge \sum_{k=1}^m \frac{\lvert G_k\rvert}{\lvert G_k\rvert} = m \]
by another application of \Cref{L:1Dec2023}. This establishes the desired upper bound on~$m$ for $C=2$ and $t=p$, thereby completing our proof that  the second supremum in~\eqref{L:1.1step2Bp:eq1} is finite.

We now turn our attention to the first  supremum in~\eqref{L:1.1step2Bp:eq1}. Assume towards a contradiction that the set $\{\lvert\theta^{-1}(n)\rvert : n\in\N\}$ is unbounded. Given a non-empty set $G\subseteq \N$, it will be convenient to introduce the notation \mbox{$G^\dagger = G\setminus\{\min G\}$}. Arguing as in the proof of \cite[Prop\-o\-si\-tion~3.11]{GL}, we can recursively construct an increasing sequence of maximal Schreier sets $G_1<G_2<\cdots$, each contained in \mbox{$M\cap[2,\infty)$}, and a sequence $(n_j)_{j\in\N}$ of natural numbers such that 
\begin{equation}\label{L:1.1step2Bp:eq4} 
\theta(i) = n_j\qquad (j\in\N,\, i\in G_j^\dagger). 
\end{equation}
We include the details of this recursion for the reader's convenience. Set $m_1 = \min M\cap [2,\infty)$. By hypothesis, we can choose a number $n_1\in\N$ such that $\lvert\theta^{-1}(n_1)\rvert > m_1$, so we can find a subset $F_1\subseteq \theta^{-1}(n_1)\setminus\{1,m_1\}$ of cardinality~$m_1-1$. Then $G_1 = \{m_1\}\cup F_1\subseteq M\cap [2,\infty)$ is a maximal Schreier set, and~\eqref{L:1.1step2Bp:eq4} is satisfied for $j=1$ because $G_1^\dagger = F_1$. 

Now assume recursively that $G_1<\cdots<G_{j-1}$ have been chosen for some $j\ge 2$. Set \mbox{$m_j = \min M\cap (\max G_{j-1},\infty)$}, choose $n_j\in\N$ such that $\lvert\theta^{-1}(n_j)\rvert\ge m_j + \lvert M\cap[1,m_j)\rvert$, and take a subset $F_j\subseteq \theta^{-1}(n_j)\cap(m_j,\infty)$ of cardinality~$m_j-1$. Then $G_j = \{m_j\}\cup F_j\subseteq M\cap [m_j,\infty)$ is a maximal Schreier set such that~\eqref{L:1.1step2Bp:eq4} is satisfied for the given value of~$j$, and $G_j>G_{j-1}$ because $\min G_j = m_j>\max G_{j-1}$. Hence the recursion continues. 

Choose an integer~$m$ such that $m>(4\lVert U\rVert)^{\frac{p}{p-1}}$ if $E=B_p$ and $m>4\lVert U\rVert^p$ if $E=S_p$. We observe that the set $\{ n_j : j\in\N\}$ is unbounded, or else we could take $n\in\N$ and $J\in[\N]$ such that $n_j = n$ for every $j\in J$, which would imply that $\bigcup_{j\in J}G_j^\dagger\subseteq\theta^{-1}(n)$, contradicting that~$\theta^{-1}(n)$ is finite, as shown in the first part of the proof. Consequently, we can find a set $K\in[\N]^{<\infty}$ such that $\lvert K\rvert = m\le\min\{ n_k : k\in K\}$, and therefore  $\{ n_k : k\in K\}\in\mathcal{S}_1$. 

For $E=B_p$, consider the vector $y = \sum_{k\in K}\lvert G_k\rvert^{-1}\sum_{j\in G_k^\dagger} e_j\in B_p$, which has norm at most $2m^{\frac1p}$ by~\eqref{estimate2baernstein} and the $1$\nobreakdash-un\-con\-di\-tion\-al\-i\-ty of the basis~$(e_j)_{j\in\N}$. We can now argue as in~\eqref{L:1.1step2Bp:eq3} to obtain
\begin{equation*}
  2^pm\lVert U\rVert^p\ge \lVert U y\rVert_{B_p}^p 
  = \biggl\lVert \sum_{k\in K}\frac{\lvert G_k\rvert-1}{\lvert G_k\rvert}e_{n_k}\biggr\rVert_{B_p}^p = \biggl(\sum_{k\in K}\frac{\lvert G_k\rvert-1}{\lvert G_k\rvert}\biggr)^p\ge\Bigl(\frac{m}{2}\Bigr)^p,
\end{equation*}
where we have used~\eqref{L:1.1step2Bp:eq4}, \Cref{L:1Dec2023} and the fact that $\lvert G_k\rvert-1\ge \lvert G_k\rvert/2$ for every $k\in\N$. Rearranging this inequality, we find $4^p\lVert U\rVert^p\ge m^{p-1}$, which contradicts our choice of~$m$. 

Again, the argument for $E=S_p$ is very similar, just using the vector \[ y = \sum_{k\in K}\frac{1}{\lvert G_k\rvert^{\frac1p}}\sum_{j\in G_k^\dagger} e_j\in S_p, \] whose norm is at most~$2^{\frac1p}$. Following the same steps as above, we obtain
\begin{equation*}
  2\lVert U\rVert^p \ge \lVert U y\rVert_{S_p}^p = \biggl\lVert \sum_{k\in K}\frac{\lvert G_k\rvert-1}{\lvert G_k\rvert^{\frac1p}}e_{n_k}\biggr\rVert_{S_p}^p\ge\sum_{k\in K}\frac{(\lvert G_k\rvert-1)^p}{\lvert G_k\rvert}\ge\frac{m}{2},
\end{equation*}
once again contradicting the choice of~$m$.
\end{proof}

\begin{proposition}\label{1.1step3Bp}
Let $M,N\in [\mathbb{N}]$, and suppose that there exists a map $\theta\colon M\rightarrow N$ for which 
\begin{equation}\label{1.1step3Bp:eq1}
    \sup\{\lvert\theta^{-1}(n)\rvert : n\in N\}<\infty\qquad\text{and}\qquad \sup \{\tau_1(\theta^{-1}(F)) : F\in\mathcal{S}_1\cap[N]^{<\infty}\}<\infty.
\end{equation}   
Then $\gl(M,N) < \infty$. 
\end{proposition}

\begin{proof}
This is a restatement of \cite[Proposition~3.12]{GL} for $\xi =1$, bearing in mind that \[ \sup\{\lvert\theta^{-1}(n)\rvert : n\in N\} = \sup\{ \tau_0(\theta^{-1}(F)) : F\in \mathcal{S}_0\cap[N]^{<\infty}\} \] 
because $\mathcal{S}_0 = \{ \{n\} : n\in\N\}\cup\{\emptyset\}$ and $\tau_0(A) = \lvert A\rvert$ for every $A\in[\N]^{<\infty}$.
\end{proof}

\begin{proof}[Proof of Theorem~{\normalfont{\ref{GLThm}}}, \ref{GLThm1iii}$\Rightarrow$\ref{GLThm1ii}] Suppose that $T\in\mathscr{B}(E_M,E_N)$ is an operator for which \[ \inf_{m\in M}\lVert Te_{m}\rVert_\infty > 0. \] Then~$T$ satisfies condition~\eqref{1.1step1BpV2:eq1} with respect to the bases~$(e_m)_{m\in M}$ and $(e_n)_{n\in N}$ for~$E_M$ and~$E_N$, respectively, so \Cref{1.1step1BpV2} shows that there is an operator $U\in \mathscr{B}(E_M,E_N)$ for which $Ue_m = e_{\theta(m)}$ for every $m\in M$, where $\theta(m)\in N$ is a suitably chosen index. Regarding~$\theta$ as a map of~$M$ into~$\N$, we can apply \Cref{L:1.1step2Bp} to deduce that both suprema in~\eqref{L:1.1step2Bp:eq1} are finite. However, they are equal to the suprema in~\eqref{1.1step3Bp:eq1} because $\theta(M)\subseteq N$, so \Cref{1.1step3Bp} implies that $\gl(M,N) < \infty$, as required.
\end{proof}

With the proof of \Cref{GLThm} complete, we state an important consequence of it that is the second main outcome of this section. 
\begin{theorem}\label{GLThmPart2}
    Let $E= B_p$ for some $1<p<\infty$ or $E=S_p$ for some $1\le p<\infty$. 
The following conditions are equivalent for $M,N\in[\N]\colon$
\begin{enumerate}[label={\normalfont{(\alph*)}}]
    \item\label{GLThm3i}  The Gasparis--Leung indices $\gl(M,N)$ and $\gl(N,M)$ are both finite. 
    \item\label{GLThm3ii} The basic sequences $(e_{m})_{m\in M}$ and $(e_{n})_{n\in N}$ are equivalent.
    \item\label{GLThm3iii} The subspaces~$E_M$ and~$E_N$ 
     are isomorphic.
\end{enumerate}    
\end{theorem}

As already indicated, we shall deduce this result from \Cref{GLThm}. However, the implication \ref{GLThm3iii}$\Rightarrow$\ref{GLThm3i} requires one additional ingredient: every isomorphic embedding of~$E_M$ into~$E$ satisfies the technical condition~\ref{GLThm1iii} of \Cref{GLThm}. 

\begin{lemma}\label{L:oldThm1.1ii} Let $E=B_p$ for some $1< p<\infty$ or $E=S_p$ for some $1\le p<\infty$, and suppose that $T\in\mathscr{B}(E_M,E)$ is an isomorphic embedding for some set $M\in[\N]$. Then $\inf_{m\in M}\lVert Te_m\rVert_\infty > 0$.
\end{lemma}

\begin{proof} Assume towards a contradiction that $\inf_{m\in M}\lVert Te_m\rVert_\infty =0$ for some $M\in[\N]$ and some isomorphic embedding $T\in\mathscr{B}(E_M,E)$. Take $\eta>0$ such that $\lVert Tx\rVert\ge \eta\lVert x\rVert$ for every $x\in E$, and set $k_{0} = 0$, $P_0=0$ and $\varepsilon_{j} = \eta/(3\cdot2^{j}+1)$ for $j \in \N$. We can then recursively choose increasing sequences $(m_{j})_{j\in\N}$ in~$M$ and $(k_{j})_{j\in\N}$ in~$\N$ such that
\begin{equation*}
\lVert Te_{m_{j}}\rVert_\infty\leq\frac{\varepsilon_{j}}{2(k_{j-1}+1)}\qquad\text{and}\qquad \lVert (I_E-P_{k_{j}})Te_{m_{j}}\rVert_E\leq \frac{\varepsilon_{j}}2\qquad (j \in \N).
\end{equation*}
This implies that for each $j\in\N$, the vector $v_{j} = (P_{k_j} - P_{k_{j-1}})Te_{m_{j}}\in E$ satisfies 
\begin{align*}
\lVert Te_{m_{j}} - v_j\rVert_E &\le \lVert (I_E - P_{k_{j}})Te_{m_{j}}\rVert_E + \lVert P_{k_{j-1}}Te_{m_{j}}\rVert_E\\ &\le \frac{\varepsilon_{j}}2 + k_{j-1}\cdot\max_{1\le n\le k_{j-1}}\lvert\langle Te_{m_{j}},e_{n}^{*}\rangle\rvert\le\frac{\varepsilon_{j}}{2} + \frac{\varepsilon_{j}}{2}= \varepsilon_{j}.
\end{align*}
In particular we have 
\[ \lVert T\rVert\ge\lVert v_{j}\rVert_E\ge  \lVert Te_{m_{j}}\rVert_E - \lVert Te_{m_{j}} - v_{j}\rVert_E\geq \eta - \varepsilon_{j} \geq \frac{6\eta}{7}, \]
so $(v_j)_{j\in\N}$ is a semi-normalized block basic sequence of $(e_n)_{n\in\N}$.
Furthermore, since 
\[ \sum_{j=1}^\infty\frac{\lVert Te_{m_{j}} - v_{j}\rVert_E}{\lVert v_{j}\rVert_E}\le \sum_{j=1}^\infty\frac{\varepsilon_j}{\eta-\varepsilon_j} = \sum_{j=1}^\infty\frac{1}{3\cdot2^{j}} = \frac13<\frac12,
\]
the Principle of Small Perturbations (see for instance \cite[Theorem~1.3.9]{AK}) implies that  $(Te_{m_{j}})_{j\in\N}$ is a basic sequence equivalent to~$(v_j)_{j\in\N}$. (Here we have used the fact that the basis constant of~$(v_j)_{j\in\N}$ is no greater than the basis constant of~$(e_n)_{n\in\N}$, which is~$1$.) 

Set $u_{j} = v_{j}/\lVert v_{j}\rVert_E$ for $j\in\N$. Being semi-normalized and unconditional, $(v_j)_{j\in\N}$ is equivalent to~$(u_j)_{j\in\N}$, and 
\[ \lVert u_j\rVert_\infty = \frac{\lVert(P_{k_j} - P_{k_{j-1}})Te_{m_{j}}\rVert_\infty}{\lVert v_{j}\rVert_E}\le\frac{\lVert Te_{m_{j}}\rVert_\infty}{6\eta/7}\le \frac{7}{12(3\cdot 2^j+1)(k_{j-1}+1)}\rightarrow 0\quad\text{as}\quad j\to\infty,  \]
so $\inf_{j\in\N}\lVert u_j\rVert_\infty =0$. Hence \Cref{P:charDsubseq} implies that $(u_j)_{j\in\N}$ admits a subsequence $(u_{j_n})_{n\in\N}$ which is equivalent to the unit vector basis $(d_n)_{n\in\N}$ for~$D$, where $D=\ell_p$ if $E=B_p$ and $D=c_0$ if $E=S_p$, as usual. 

In conclusion, we have shown that $(d_n)_{n\in\N}$ is equivalent to $(u_{j_n})_{n\in\N}$, which is equivalent to $(v_{j_n})_{n\in\N}$, which is equivalent to $(Te_{m_{j_n}})_{n\in\N}$, and therefore $(e_{m_{j_n}})_{n\in\N}$ is equivalent to $(d_n)_{n\in\N}$ because~$T$ is an isomorphic embedding. However, this is absurd: no subsequence of $(e_n)_{n\in\N}$ is dominated by $(d_n)_{n\in\N}$, as is easy to see (or alternatively this is a very special case of \Cref{P:charDsubseq}).
\end{proof}

\begin{corollary}\label{C:9aug2024}
    Let $E=B_p$ for some $1< p<\infty$ or $E=S_p$ for some $1\le p<\infty$, and suppose that~$E_M$ embeds isomorphically into~$E_N$ for some sets $M,N\in[\N]$. Then $\gl(M,N)<\infty$.
\end{corollary}

\begin{proof}
   Take an isomorphic embedding $T\in\mathscr{B}(E_M,E_N)$, and let $J\colon E_N\to E$ denote the inclusion map. Then we have $0 < \inf_{m\in M}\lVert JTe_m\rVert_\infty = \inf_{m\in M}\lVert Te_m\rVert_\infty$ by \Cref{L:oldThm1.1ii}, so the implication \ref{GLThm1iii}$\Rightarrow$\ref{GLThm1ii} in \Cref{GLThm} shows that $\gl(M,N)<\infty$.
\end{proof}

\begin{proof}[Proof of Theorem~{\normalfont{\ref{GLThmPart2}}}]
The equivalence of conditions~\ref{GLThm1ii} and~\ref{GLThm1i} in \Cref{GLThm} implies that conditions~\ref{GLThm3i} and~\ref{GLThm3ii} are also equivalent in \Cref{GLThmPart2}. The implication \ref{GLThm3ii}$\Rightarrow$\ref{GLThm3iii} is clear, and finally \Cref{C:9aug2024} shows that \ref{GLThm3iii} implies~\ref{GLThm3i}. 
\end{proof}

We conclude this section with two applications of \Cref{C:9aug2024}, both establishing counter\-parts for the Baernstein and Schreier spaces of results of Gasparis and Leung concerning the structure of the complemented subspaces of the higher-order Schreier spaces. 

\begin{definition} \begin{enumerate}[label={\normalfont{(\roman*)}}]   
\item A Banach space~$X$ is \emph{primary} if the kernel or the range of~$P$ is isomorphic to~$X$ for every idempotent operator $P\in\mathscr{B}(X)$.    
\item Two Banach spaces~$X$ and~$Y$ are \emph{incomparable} if no subspace of~$X$ is isomorphic to~$Y$ and  no subspace of~$Y$ is isomorphic to~$X$. 
\end{enumerate}
\end{definition}

\begin{proposition}\label{GL:subspaces}
Let $E = B_p$ for some $1<p<\infty$ or $E=S_p$ for some $1\le p<\infty$. Then:
\begin{enumerate}[label={\normalfont{(\roman*)}}]   
\item\label{GL:subspaces:1} The subspace~$E_N$ fails to be primary for every $N\in[\N]$.
\item\label{GL:subspaces:2} There is a subset~$\mathcal{A}$ of~$[\N]$ of cardinality~$\mathfrak{c}$ such that~$E_{L}$ and~$E_{M}$ are incomparable whenever $L,M\in\mathcal{A}$ are distinct.
\end{enumerate}
\end{proposition}

\begin{proof} We follow the approach Gasparis and Leung took in their proofs of \cite[Corollary~3.15 and Theorem~1.3]{GL}, respectively. The common starting point is that for any set $N\in[\N]$, we can equip~$[N]$ with the topology of pointwise convergence obtained by identifying the elements of~$[N]$ with their indicator functions; this turns~$[N]$ into a Polish space. 

\ref{GL:subspaces:1}. The set     
    \[ \mathcal{F} = \{ (L,M)\in[N]\times[N] : L\cup M = N,\, L\cap M=\emptyset\} \]
    is closed with respect to the product topology on~$[N]\times[N]$, and 
    \[ \mathcal{G} = \{ (L,M)\in\mathcal{F} : \gl(N,L) = \gl(N,M) = \infty \} \]
    is a dense $G_\delta$-subset.  In particular~$\mathcal{G}$ is non-empty, so we can take $(L,M)\in\mathcal{G}$. The fact that $N=L\cup M$ and $L\cap M=\emptyset$ implies that $E_N = E_L\oplus E_M$, but~$E_N$ is neither isomorphic to~$E_L$ nor~$E_M$ by \Cref{C:9aug2024} because  $\gl(N,L) = \gl(N,M) = \infty$. (In fact, $E_L$ and~$E_M$ do not even contain subspaces which are isomorphic to~$E_N$.) This proves that~$E_N$ is not primary. 

    \ref{GL:subspaces:2}. By \cite[Lemma~3.5]{GL}, 
    \[ \mathcal{D} = \{ (L,M)\in[\N]\times[\N] : \gl(L,M) = \gl(M,L) = \infty \} \]
     is a dense $G_\delta$-subset of~$[\N]\times[\N]$. Therefore, applying \cite[Proposition~3.6]{GL}, we can find a sub\-set~$\mathcal{A}$ of~$[\N]$ which is homeomorphic to the Cantor set and satisfies $(L,M)\in\mathcal{D}$ whenever $L,M\in\mathcal{A}$ are distinct. In particular~$\mathcal{A}$ has cardinality~$\mathfrak{c}$, and  \Cref{C:9aug2024} shows that~$E_L$ and~$E_M$ are incomparable for distinct $L,M\in\mathcal{A}$ because $\gl(L,M) = \gl(M,L) = \infty$.
 \end{proof}

\section{Spatial ideals of operators on the Baernstein and Schreier spaces}\label{S:5}
\noindent Let $X$ be a Banach space. We write $\langle T\rangle$ for the (algebraic, two-sided) ideal of~$\mathscr{B}(X)$
generated by an operator $T\in\mathscr{B}(X)$, that is,
\begin{equation}\label{eq:classicalideal}
    \langle T\rangle = \biggl\{\sum_{j=1}^k U_j TV_j : k\in\N,\, U_1,\ldots,U_k,V_1,\ldots,V_k\in\mathscr{B}(X)\biggr\}.
\end{equation}  
Since $\mathscr{B}(X)$ is a unital Banach algebra, the ideal~$\langle
T\rangle$ is proper if and only if its norm-closure $\overline{\langle T\rangle}$ is. 
Suppose that~$X$ has an unconditional basis. 
Following~\cite{BKL}, we call the closed ideals of the form $\overline{\langle P_M\rangle}$ for some non-empty subset~$M$ of~$\N$ \emph{spatial}, where~$P_M$ denotes the basis projection, as usual.  

The main aim of this section is to prove the following proposition, which is an extended counter\-part of \cite[Proposition~4.12]{BKL} for the Baernstein and Schreier spaces. The key difference is the addition of 
a new quantitative condition, \ref{3blemmaiv}, that will play an essential role in the proofs of parts~\ref{itemsmall} and~\ref{itemlarge} of \Cref{thmsmalllargeideals} in the next section.

\begin{proposition}\label{lemmaonpage3b}
  Let $E=B_{p}$ for some $1<p<\infty$ or $E=S_{p}$ for some $1\le p<\infty$. The following conditions are equivalent for every pair of sets $M\subseteq \N$ and $N \in [\N]\colon$
\begin{enumerate}[label={\normalfont{(\alph*)}}]
\item\label{3blemmai}$P_{M} \in \overline{\langle P_{N}\rangle}$,
\item\label{3blemmaii} $\langle P_{M} \rangle \subseteq \langle P_{N} \rangle$,
\item\label{3blemmaiii} $\langle P_{N} \rangle = \langle P_{M \cup N} \rangle$,
\item\label{3blemmaiv} $\operatorname{dist}(P_{M},\langle P_{N} \rangle) < 1$,
\item\label{3blemmav} $\gl(M \cup N,N) < \infty$,
\item\label{3blemmavi} $E_{N}$ contains a subspace which is isomorphic to $E_{M}$,
\item\label{3blemmavii} $E_{N}$ contains a complemented subspace which is isomorphic to $E_{M}$,
\item\label{3blemmaviii} $E_{N}$ is isomorphic to $E_{M \cup N}$,
\item\label{3blemmaix} The basic sequences $(e_{n})_{n \in M \cup N}$ and $(e_{n})_{n \in N}$ are equivalent.
\end{enumerate}
\end{proposition}

We require two lemmas in the proof of this proposition. The first is a variant of the Neumann series, showing that every idempotent element which is close to an ideal of a Banach algebra must in fact belong to the ideal.  

\begin{lemma}\label{lemmaonpage3a}
    Let $\mathscr{I}$ be an ideal of a Banach algebra $\mathscr{A}$, and take a non-zero idempotent $p \in \mathscr{A}$. Then $p \in \mathscr{I}$ if (and only if) $\operatorname{dist}(p,\mathscr{I}) < \lv p \rv^{-2}$.
\end{lemma}

\begin{proof} The implication $\Rightarrow$ is obvious. Conversely, suppose that $\lv p - a \rv < \lv p \rv^{-2}$ for some $a \in \mathscr{I}$. Then $\lv p - pap \rv < 1$, so the series $\sum_{n=1}^{\infty}(p-pap)^{n}$ converges absolutely. Set $b = p + \sum_{n=1}^{\infty}(p-pap)^{n} \in \mathscr{A}$ and observe that 
    \begin{align*}
      \mathscr{I}\ni bpap &= \biggl(p + \sum_{n=1}^{\infty}(p-pap)^{n}\biggr)\bigl(p - (p-pap)\bigr)\\ &= p - p(p-pap) + \sum_{n=1}^{\infty}(p-pap)^{n}p - \sum_{n=2}^{\infty}(p-pap)^{n} = p.\qedhere
    \end{align*}
\end{proof}

The second lemma is the counterpart of \cite[Proposition~4.6]{BKL}. It will enable us to connect the first four conditions of \Cref{lemmaonpage3b} concerning ideals with the last four (or five) concerning subspaces. 

\begin{lemma}\label{cartesiansquare}
Let $E= B_{p}$ for some $1<p<\infty$ or $E=S_{p}$ for some $1\le p<\infty$. Then 
\[ E_{M}\cong E_{M} \oplus E_{M}\qquad (M \in [\N]). \]
\end{lemma}

\begin{proof}
    We take the same approach as in the proof of \cite[Proposition~4.6]{BKL}. Set $M' = 2M-1$ and $M'' = 2M$. Since these sets are disjoint, we have
    \[ E_{M'\cup M''}  = E_{M'}\oplus E_{M''} \]
    by unconditionality, so it will suffice to show that each of these spaces is isomorphic to~$E_M$, which in turn will follow from \Cref{GLThmPart2} provided that the appropriate Gasparis--Leung indices are finite. First, we have $\gl(M,M'\cup M'')\le 3$ and $\gl(M'\cup M'', M)\le 2$ by \cite[Lemma~4.11]{BKL}, so $E_{M'\cup M''}\cong E_M$. Second, $M''$ is a spread of~$M$ and of~$M'$, so $\gl(M'',M)=1=\gl(M'',M')$.  Third, we claim that
    \begin{equation}\label{cartesiansquare:eq1}
        \gl(M,M'')\le 2\qquad\text{and}\qquad \gl(M',M'')\le 2. 
    \end{equation} 
    Indeed, suppose that $M''(J)\in\mathcal{S}_1$ for some non-empty $J\in[\N]^{<\infty}$. Then $2m_{j_1}\ge k$, where we have written $J = \{j_1<\cdots<j_k\}$ and $M = \{m_1<m_2<\cdots\}$. This implies that we can partition~$J$ into two subsets, $J_1$ and $J_2$, each having at most~$m_{j_1}$ elements, and therefore $M(J_1),M(J_2),M'(J_1),M'(J_2)\in\mathcal{S}_1$.  Hence we have $\tau_1(M(J))\le 2$ and $\tau_1(M'(J))\le 2$ because $M(J) = M(J_1)\cup M(J_2)$ and $M'(J) = M'(J_1)\cup M'(J_2)$. This proves~\eqref{cartesiansquare:eq1}, and consequently $E_M\cong E_{M''}\cong E_{M'}$. 
\end{proof}

\begin{proof}[Proof of Proposition~{\normalfont{\ref{lemmaonpage3b}}}]
    \Cref{cartesiansquare} implies that $E_{N} \cong E_{N}\oplus E_{N}$ and $E_{M \cup N} \cong E_{M \cup N} \oplus E_{M \cup N}$ because~$N$ is infinite. Hence conditions~\ref{3blemmai}, \ref{3blemmaii}, \ref{3blemmaiii}, \ref{3blemmavii} and~\ref{3blemmaviii} are equivalent by \cite[Lemma~2.3 and Corollary~2.5]{BKL}. 
    
    We have $\gl(N,M \cup N) = 1$ because~$N$ is a spread of~$M\cup N$, so conditions~\ref{3blemmav}, \ref{3blemmaviii} and~\ref{3blemmaix} are equivalent by \Cref{GLThmPart2}. 
    
    The implications \ref{3blemmai}$\Rightarrow$\ref{3blemmaiv} and \ref{3blemmavii}$\Rightarrow$\ref{3blemmavi} are trivial, while \Cref{lemmaonpage3a} shows that~\ref{3blemmaiv} implies~\ref{3blemmaii}. We complete the proof by showing that \ref{3blemmavi} implies~\ref{3blemmav}. Suppose that $E_{M}$ embeds isomorphically into~$E_{N}$. Then $E_{M\cup N} = E_{M} \oplus E_{N\setminus M}$ embeds isomorphically into $E_{N} \oplus E_{N} \cong E_{N}$, so $\gl(M \cup N,N) < \infty$ by \Cref{C:9aug2024}.        
  \end{proof}

\Cref{lemmaonpage3b} enables us to establish a counterpart for the Baernstein and Schreier spaces of the main result of~\cite{BKL}. This requires one additional piece of terminology. Let~$X$ be a Banach space with an unconditional basis. The ideal~$\mathscr{K}(X)$ of compact operators is al\-ways spatial because $\overline{\langle P_M\rangle}=\mathscr{K}(X)$ if (and only if) $M\in[\N]^{<\infty}\setminus\{\emptyset\}$. Following~\cite{BKL}, we call a spatial ideal~$\mathscr{I}$ \emph{non-trivial} if $\mathscr{K}(X)\subsetneq\mathscr{I}\subsetneq\mathscr{B}(X)$. 

\begin{theorem}\label{manychainsofideals}
 Let $E = B_p$ for some $1<p<\infty$ or $E=S_p$ for some $1\le p<\infty$. Then:
\begin{enumerate}[label={\normalfont{(\roman*)}}]
    \item\label{manychainsofideals1} The family of non-trivial spatial ideals of $\mathscr{B}(E)$ is non-empty and has no minimal or maximal elements.
    \item\label{manychainsofideals2} Let $\mathscr{I} \subsetneq \mathscr{J}$ be spatial ideals of $\mathscr{B}(E)$. Then there is a family $\{\Gamma_{L}:L \in \Delta\}$ such that:
\begin{enumerate}[label={\normalfont{(\arabic*)}}]
	\item\label{manychainsofideals2.1} the index set $\Delta$ has the cardinality of the continuum;
	\item\label{manychainsofideals2.2} for each $L \in \Delta, \Gamma_{L}$ is an uncountable chain of spatial ideals of $\mathscr{B}(E)$ such that
	\[	\mathscr{I} \subsetneq \mathscr{L} \subsetneq \mathscr{J} \qquad (\mathscr{L} \in \Gamma_{L}),
	\] 
	and $\bigcup \Gamma_{L}$ is a closed ideal that is not spatial;
	\item\label{manychainsofideals2.3} $\overline{\mathscr{L} + \mathscr{M}} = \mathscr{J}$ whenever $\mathscr{L} \in \Gamma_{L}$ and $\mathscr{M} \in \Gamma_{M}$ for distinct $L,M \in \Delta$.
    \end{enumerate}
    \item\label{manychainsofideals3} The Banach algebra $\mathscr{B}(E)$ contains at least continuum many maximal ideals.
    \item\label{manychainsofideals4} The ideal 
    \[ \bigcap \{\mathscr{I} : \mathscr{I}\ \text{is a non-trivial spatial ideal of}\ \mathscr{B}(E)\} \] is not contained in the ideal of strictly singular operators on~$E$.
\end{enumerate}
\end{theorem}

\begin{proof} 
Clauses~\ref{manychainsofideals1}--\ref{manychainsofideals3} are the counterparts for~$E$ of \cite[Theorem~1.1]{BKL} and can be proved in exactly the same way; see \cite[pages 10--11]{BKL}. This requires that we establish the counterpart of \cite[Lemma~2.8]{BKL} for~$E$, which we can do by copying the proof given in \cite[pages 21--24]{BKL} for $n=1$, just referring to \Cref{lemmaonpage3b} instead of \cite[Proposition~4.12]{BKL} throughout.

\ref{manychainsofideals4} is the counterpart of \cite[Theorem~1.2(ii), Equation~(1.1)]{BKL}, and the proof is similar. Indeed, let $D=\ell_p$ if $E=B_p$  and $D=c_0$ if $E=S_p$, and take a projection $Q\in\mathscr{B}(E)$ whose range is isomorphic to~$D$. Then~$Q$ is not strictly singular, but \Cref{lpsaturation} implies that $Q\in\langle P_N\rangle$ for every $N\in[\N]$, and therefore~$Q$ belongs to every non-trivial spatial ideal of~$\mathscr{B}(E)$. 
\end{proof}

\section{Finding $2^{\mathfrak{c}}$ many closed ideals of operators: the proofs of \Cref{thmsmalllargeideals}\ref{itemsmall} and~\ref{itemlarge}}\label{S:6}
\noindent In this section we combine our previous results about the  Gasparis--Leung index and ideals generated by basis projections to prove the remaining two parts of \Cref{thmsmalllargeideals}; that is, taking $E=B_p$ for $1<p<\infty$ or $E=S_p$ for $1\le p <\infty$, as usual, we shall show that~$\mathscr{B}(E)$ contains~$2^{\mathfrak{c}}$ many closed ideals that lie between the ideals of compact and strictly singular operators, as well as~$2^{\mathfrak{c}}$ many closed ideals that are ``large'' in the sense that they contain projections of infinite rank. Note that~$\mathscr{B}(E)$ cannot contain more than~$2^{\mathfrak{c}}$ many closed ideals because~$E$ is separable. 

Both results rely on a general theorem of Freeman, Schlumprecht and Zs\'{a}k~\cite[Proposition~1]{FSZ} that extracts the key idea of the argument that Johnson and Schechtman~\cite{JS} used to show that~$\mathscr{B}(L_p[0,1])$ contains~$2^{\mathfrak{c}}$ many closed ideals for every $p\in(1,2)\cup(2,\infty)$.
Before we can state the said theorem of Freeman--Schlumprecht--Zs\'{a}k precisely, we require two additional pieces of terminology. The first generalizes the classical notion of a (closed) ideal of a (Banach) algebra to the space of operators between two distinct Banach spaces~$X$ and~$Y$: 
a \emph{(closed) ideal} of~$\mathscr{B}(X,Y)$ is a (norm-closed) subspace~$\mathscr{J}$ of~$\mathscr{B}(X,Y)$ such that $UTV \in \mathscr{J}$ whenever $V\in \mathscr{B}(X)$, $T\in\mathscr{J}$ and $U\in \mathscr{B}(Y)$. Extending~\eqref{eq:classicalideal}, for a subset~$\mathscr{T}$ of~$\mathscr{B}(X,Y)$, we write~$\langle\mathscr{T}\rangle$ for the ideal of~$\mathscr{B}(X,Y)$ it generates. 

The reason this generalization is useful for our purposes is that in the case where~$X$ contains a complemented subspace isomorphic to~$Y$, the map 
\begin{equation}\label{idealsetup}
    \mathscr{J} \mapsto \overline{\biggl\{\sum_{j=1}^{n}U_{j}T_{j} : n \in \N,\, U_{1},\ldots,U_{n}\in \mathscr{B}(Y,X),\,  T_{1},\ldots,T_{n} \in \mathscr{J}\biggr\}}
\end{equation}
is an injection from the lattice of closed ideals of $\mathscr{B}(X,Y)$ into the lattice of closed ideals of~$\mathscr{B}(X)$. (This is a special case of an observation stated above \cite[Proposition~2]{FSZ}, and is also easy to verify directly.) Hence, to show that~$\mathscr{B}(X)$ contains~$2^{\mathfrak{c}}$ many closed ideals, it suffices to find a complemented subspace~$Y$ of~$X$ for which~$\mathscr{B}(X,Y)$ contains~$2^{\mathfrak{c}}$ many closed ideals. 

The second notion that we require is that of a \emph{$1$-un\-con\-di\-tional finite-di\-men\-sion\-al decomposition,} or \emph{$1$-UFDD} for short, of a Banach space~$X$; that is, a sequence $(X_n)_{n\in\N}$ of finite-di\-men\-sional subspaces of~$X$ such that every $x\in X$ has a unique decomposition of the form $x=\sum_{n=1}^\infty x_n$, where $x_n\in X_n$ for every $n\in\N$, and the series $\sum_{n=1}^\infty\sigma_n x_n$ converges with $\lVert\sum_{n=1}^\infty\sigma_n x_n\rVert \le \lVert x\rVert$ for every sequence $(\sigma_n)_{n\in\N}\in\{\pm 1\}^\N$. It follows that for every (non-empty) subset~$N$ of~$\N$, we can define a projection $Q_N\in\mathscr{B}(X)$ of norm~$1$ by $Q_Nx = \sum_{n\in N}x_n$. 

In fact, we shall only consider $1$-UFDDs of a very simple kind. Let~$X$ be a Banach space with a $1$-un\-condi\-tional basis $(x_n)_{n\in\N}$, and take a partition $J_1<J_2<\cdots$ of~$\N$ into finite, successive intervals. Then the sequence of finite-dimensional subspaces given by
\begin{equation}\label{eq:UFDD} X_n = \operatorname{span}\{ x_j : j\in J_n\}\qquad (n\in\N) \end{equation} 
is a $1$-UFDD for~$X$. 
For later reference, we observe that in this case the projection~$Q_N$, for $N\subseteq\N$, defined above is equal to the basis projection~$P_{L_N}$ induced by the set $L_N = \bigcup_{n\in N}J_n$.

\begin{theorem}[Freeman, Schlumprecht and Zs\'{a}k]\label{thm:FSZ} Let $\mathcal{A}\subset [\N]$  be an almost disjoint family of cardinality~$\mathfrak{c}$, let~$X$ and~$Y$ be Banach spaces with $1$-UFDDs $(X_n)_{n\in\N}$ and $(Y_n)_{n\in\N}$, respectively, and suppose that $T\in\mathscr{B}(X,Y)$ is an operator which satisfies 
\begin{enumerate}[label={\normalfont{(\roman*)}}]
\item\label{thm:FSZ:i} $T[X_n]\subseteq Y_n$ for every $n\in\N;$
\item\label{thm:FSZ:ii} $\inf\bigl\{\operatorname{dist}(TQ_M,\langle TQ_N\rangle) : M,N\in[\N],\, \lvert M\setminus N\rvert=\infty\bigr\} >0$. 
\end{enumerate}
Then the map 
\begin{equation}\label{thm:FSZ:eq1}
    \mathcal{N}\mapsto \overline{\langle TQ_N : N\in\mathcal{N}\rangle}
\end{equation}  
defines an order-preserving injection from the power set of~$\mathcal{A}$ into the lattice of closed ideals of~$\mathscr{B}(X,Y)$. 
\end{theorem}

To enable us to apply this theorem 
to the Baernstein and Schreier spaces, we present a variant not involving dyadic trees of the key construction that Manoussakis and Pelczar-Barwacz used in their proof of \cite[Lemma~4.3]{MPB}. 
\begin{construction}\label{ConstructionMPB}
Set $F_1=\emptyset$. By recursion, we can partition $\N$ into finite, successive intervals $G_{1} < F_{2} < G_{2} < F_{3} < G_{3} < \cdots$ with the following properties:
\begin{enumerate}[label={\normalfont{(\roman*)}}]
    \item\label{MPBsetG} $G_{n}$ is the union of $n$ successive maximal Schreier sets  for each $n \in \N$ (so in particular $\tau_{1}(G_{n}) = n$),
    \item\label{MPBsetH} $\lvert F_{n}\rvert = \sum_{m=1}^{n-1}(\lvert F_{m}\rvert + \lvert G_{m}\rvert)$ for each $n\ge 2$.
\end{enumerate}
For brevity, we introduce the notation $J_{n} = F_{n}\cup G_{n}$ for $n\in\N$ and set
\begin{equation}\label{eq:defnLM}
    L_{N} = \bigcup_{n\in N}J_{n}\qquad (N\subseteq\N).
\end{equation} 
We observe for later reference that we can rewrite property~\ref{MPBsetH} as $\lvert F_{n}\rvert = \sum_{m=1}^{n-1}\lvert J_{m}\rvert$.
\end{construction}
\begin{lemma}\label{lemmaon1b}
    For $M,N \in [\N]$, the set~$M \setminus N$ is bounded above by~$\gl(L_{M},L_{N})\in\N\cup\{\infty\}$, where $L_{M},L_{N}\in[\N]$ are defined by~\eqref{eq:defnLM}.
\end{lemma}
\begin{proof} Take $m\in M \setminus N$. We seek to prove that $m\le \gl(L_M,L_N)$, which by the definition~\eqref{Eq:GLindex} of the Gasparis--Leung index means that we must find a set $K \in [\N]^{<\infty}$ such that
    \begin{gather}\label{equation1}
        \tau_{1}(L_{M}(K)) \geq m\qquad \text{and}\qquad L_{N}(K) \in \Ss.
    \end{gather}
    The case $m=1$ is trivial, so we may suppose that $m\geq 2$. We have $J_{m}\subseteq L_{M}$ because $m\in M$, so we can find $H\in [\N]^{<\infty}$ such that $J_{m} = L_{M}(H)$. Note that~$H$ is an interval because $J_{m}$ is, and the definition of~$J_m$ implies that we can split $H$ in two subintervals $H_{1} < H_{2}$ such that $L_{M}(H_{1})  = F_{m}$ and $L_{M}(H_{2}) = G_{m}$. We claim that $K = H_{2}$ satisfies~\eqref{equation1}.
    
    The first part is immediate because $\tau_{1}(L_{M}(H_2)) = \tau_{1}(G_{m}) = m$, so it only remains to show that $L_{N}(H_2) \in \Ss$; that is, $|L_{N}(H_2)| \leq \min(L_{N}(H_2))$.  Set $k = \min H_2$ and observe that $\min(L_{N}(H_2))$ is the $k^{\text{th}}$ element of the set 
    \begin{gather*}
        L_{N} = \bigcup_{n\in N}J_{n} = \Bigl(\bigcup_{n\in N\cap [1,m)}J_{n}\Bigr)\cup\Bigl(\bigcup_{n\in N\cap (m,\infty)}J_{n}\Bigr)
    \end{gather*}
    because $m\notin N$. We have
    \begin{gather*}
        \Bigl|\bigcup_{n\in N \cap [1,m)}J_{n}\Bigr| = \sum_{n\in N \cap [1,m)}|J_{n}| \leq \sum_{j=1}^{m-1}|J_{j}| = |F_{m}| = |L_{M}(H_{1})| = |H_{1}| \leq \max H_{1} < k,
    \end{gather*}
    so the $k^{\text{th}}$ element of $L_{N}$ must belong to the set $\bigcup_{n\in N \cap (m,\infty)}J_{n}$. Hence
    \begin{align*}
        \min(L_{N}(H_2)) &\geq \min\Bigl(\bigcup_{n\in N\cap (m,\infty)}\!\!J_{n}\Bigr)\geq \min J_{m+1} > \max J_{m}\\ &= \max G_m\geq |G_{m}| = |L_{M}(H_2)| = |H_2|  = |L_{N}(H_2)|, 
    \end{align*}
    and the conclusion follows.
\end{proof}

\begin{corollary}\label{corollaryonpageB}
The following conditions are equivalent for $M,N\in [\N]\colon$
    \begin{enumerate}[label={\normalfont{(\alph*)}}]
        \item\label{pageBa} $\gl(L_{M},L_{N}) < \infty$,
        \item\label{pageBb} $\gl(L_{M\cup N},L_{N}) < \infty$,
        \item\label{pageBc} $|M \setminus N| < \infty$.
    \end{enumerate} 
\end{corollary}

\begin{proof}
  \Cref{lemmaon1b} shows that the set~$M \setminus N$ is bounded above by $\gl(L_{M},L_{N})$, so~\ref{pageBa} implies~\ref{pageBc}.

  \ref{pageBc}$\Rightarrow$\ref{pageBb}. Suppose that $M \setminus N$ is finite. Then $L_{M \cup N} \setminus L_{N} = L_{M \setminus N}$ is finite, too, so $\gl(L_{M\cup N},L_{N}) < \infty$.
  
  \ref{pageBb}$\Rightarrow$\ref{pageBa}. This is a consequence of the fact that $\Ss$ is closed under spreading.
\end{proof}

\begin{proof}[Proof of Theorem {\normalfont{\ref{thmsmalllargeideals}\ref{itemlarge}}}]
    We shall apply \Cref{thm:FSZ} with $X=Y=E$, $T = I_{E}$ and the $1$-UFDDs given by $X_n = Y_n = \operatorname{span}\{ e_j : j\in J_n\}$ for every $n\in\N$, where $J_1<J_2<\cdots$ are the intervals defined in \Cref{ConstructionMPB}. These choices ensure that condition~\ref{thm:FSZ:i} of \Cref{thm:FSZ} is trivially satisfied. To verify condition~\ref{thm:FSZ:ii}, we recall that for $N\subseteq\N$, the projection~$Q_N$ associated with the chosen $1$-UFDDs is the basis projection~$P_{L_N}$. Taking $M,N\in[\N]$ with $\lvert M\setminus N\rvert = \infty$, we have $\gl(L_M\cup L_N,L_N) = \gl(L_{M\cup N}, L_N)=\infty$ by \Cref{corollaryonpageB}, so \Cref{lemmaonpage3b} implies that
    $\operatorname{dist}(Q_M,\langle Q_N\rangle) = \operatorname{dist}(P_{L_{M}},\langle P_{L_{N}} \rangle) = 1$.
    Hence \Cref{thm:FSZ} shows that \[ \bigl\{\overline{\langle P_{L_N} : N\in\mathcal{N}\rangle} : \mathcal{N}\subseteq\mathcal{A}\bigr\} \] is a collection of $2^\mathfrak{c}$ many distinct closed ideals of $\mathscr{B}(E)$ for any almost disjoint family $\mathcal{A} \subset [\N]$ of cardinality~$\mathfrak{c}$. 

    The set $\{ UV : U\in\mathscr{B}(D,E),\, V\in\mathscr{B}(E,D)\}$ is closed under addition and therefore an ideal of~$\mathscr{B}(E)$ because $D\cong D\oplus D$, where we recall that $D=\ell_p$ if $E=B_p$ and $D=c_0$ if $E=S_p$, as usual. An easy standard argument shows that this ideal is equal to~$\langle Q\rangle$ for any projection $Q\in\mathscr{B}(E)$ whose range is isomorphic to~$D$. As we saw in the proof of \Cref{manychainsofideals}\ref{manychainsofideals4}, $Q\in\langle P_N\rangle$ for every $N\in[\N]$, so $\langle Q\rangle\subseteq \overline{\langle P_{L_N} : N\in\mathcal{N}\rangle}$ for every non-empty subset~$\mathcal{N}$ of~$\mathcal{A}$.  
\end{proof}

The proof of \Cref{thmsmalllargeideals}\ref{itemsmall} follows a similar path, but some aspects require additional work. We begin with a standard characterization of strictly singular operators defined on a Banach space with a basis, and include a short proof for completeness. 

\begin{lemma}\label{L:SSblocks}
 Let $T\in\mathscr{B}(X,Y)$ be an operator between Banach spaces~$X$ and~$Y$, and suppose that~$X$ has a basis. Then $T$ is strictly singular if (and only if) the restriction of~$T$ to any block subspace of~$X$ fails to be an isomorphic embedding. 
\end{lemma}

\begin{proof} The implication $\Rightarrow$ is trivial because block subspaces are in\-finite-di\-men\-sional. 

Conversely, suppose that~$T$ fails to be strictly singular, so that its restriction to some closed, infinite-dimensional subspace~$Z$ of~$X$ is bounded below by some number $\eta>0$. We use the same notation and approach as in the first part of the proof of \Cref{C2afterC1}; that is, $(x_{n})_{n\in\N}$ denotes the basis of~$X$,  $K$~is the basis constant, $P_{n}$ is the $n^{\text{th}}$ basis projection for $n\in\N$, and we set $m_{0} = 0$ and $P_{0} = 0$. By recursion, we choose natural numbers $m_{1} < m_{2} < \cdots$ and unit vectors $z_{n} \in Z\cap \ker P_{m_{n-1}}$ such that \[ \lv z_{n} - w_{n}\rv\leq \varepsilon_n,\quad \text{where}\quad w_{n} = P_{m_{n}}z_{n}\quad\text{and}\quad \varepsilon_n = \frac{\eta}{2^{n+2}K(\eta+\lVert T\rVert)+\eta}\qquad (n \in \N). \] 
Then $(w_{n})_{n \in \N}$ is a block basic sequence of $(x_{n})_{n \in \N}$ because $\lv w_{n} \rv\ge 1-\varepsilon_n>0$ for every $n\in\N$. 

We shall now complete the proof by showing that the restriction of~$T$ to the block subspace spanned by $(w_{n})_{n \in \N}$ is bounded below by $\eta/2$. Take  a unit vector $w = \sum_{n=1}^N\alpha_nw_n$ for some $N\in\N$ and $\alpha_1,\ldots,\alpha_N\in\K$, and set $z = \sum_{n=1}^N\alpha_nz_n\in Z$. We have 
\begin{align}\label{L:SSblocks:eq1}
    \lVert Tw\rVert &\ge \lVert Tz\rVert - \lVert T(z-w)\rVert\ge \eta\lVert z\rVert - \lVert T\rVert \lVert z-w\rVert\\ &\ge \eta(\lVert w\rVert - \lVert z-w\rVert) - \lVert T\rVert \lVert z-w\rVert = \eta - (\eta +\lVert T\rVert) \lVert z-w\rVert.\notag
\end{align}
To find an upper bound on $\lVert z-w\rVert$, we observe that $\alpha_nw_n = (P_{m_n}-P_{m_{n-1}})w$, so 
\[ \lvert\alpha_n\rvert\le \frac{2K}{\lVert w_n\rVert}\le \frac{2K}{1-\varepsilon_n}\qquad (1\le n\le N), \]
and therefore 
\begin{align*} \lVert z-w\rVert &\le \sum_{n=1}^N\lvert\alpha_n\rvert\lVert z_n-w_n\rVert\le \sum_{n=1}^N\frac{2K\varepsilon_n}{1-\varepsilon_n} = \sum_{n=1}^N\frac{\eta}{(\eta+\lVert T\rVert)2^{n+1}}\le \frac{\eta}{2(\eta+\lVert T\rVert)}, 
\end{align*}
where the equality in the middle follows from the choice of~$\varepsilon_n$. Substituting this estimate into~\eqref{L:SSblocks:eq1}, we conclude that  $\lVert Tw\rVert\ge\eta/2$, which establishes the result.
\end{proof}

\begin{proposition}\label{P:formalinclusionmap}
Let $(E,D) = (B_{p},\ell_{p})$ for some $1 < p < \infty$ or $(E,D) = (S_{p},c_{0})$ for some $1 \leq p < \infty$, and let $(e_n)_{n\in\N}$ and $(d_n)_{n\in\N}$ denote the unit vector bases for~$E$ and~$D$, respectively.   Then the formal inclusion map given by $\iota\colon e_n\mapsto d_n$ for $n\in\N$ 
extends to a bounded linear injection $\iota\colon E\to D$ of norm~$1$. Furthermore, $\iota$ is strictly singular, but not compact. 
\end{proposition}

\begin{proof}  It is obvious that the formal inclusion map $\iota\colon S_p\to c_0$ is a bounded linear injection of norm~$1$, while the same conclusion for $\iota\colon B_p\to\ell_p$ is an easy consequence of the definition of the norm on~$B_{p}$, or alternatively it follows by applying \Cref{Sigmaoperator} to the chain $\mathcal{C} = \bigl\{\{n\} : n\in\N\bigr\}$. The non-com\-pact\-ness of~$\iota$ is witnessed by its action on the unit vector basis in both cases, so it only remains to verify that~$\iota$ is strictly singular. 

\Cref{L:SSblocks} implies that it suffices to show that the restriction of~$\iota$ to the closed subspace spanned by a block basic sequence $(w_{n})_{n\in\N}$ of~$(e_n)_{n\in\N}$ is not an isomorphic embedding. By \Cref{infinitynull}, $(w_n)_{n\in\N}$ admits a normalized block basic sequence $(u_{n})_{n\in\N}$ for which $\lv u_{n} \rv_{\infty} \rightarrow 0$ as $n\to\infty$. This completes the proof for the Schreier space~$S_p$ because $\lVert \iota(u_n)\rVert = \lv u_{n} \rv_{\infty}$ in this case. The argument for the Baernstein space~$B_p$ is more subtle, relying on an inequality due to Jameson that we shall establish in \Cref{App:GJOJ} below; it involves a constant $K_p>0$ which depends only on~$p$. Using the variant of Jameson's inequality stated in the last line of \Cref{ThmgGJOJp}, we obtain 
\[  \lv \iota(u_{n}) \rv^p = \lv u_n\rv_{\ell_p}^p \leq K_{p}\lv u_{n} \rv_{\infty}^{p-1}\lv u_{n} \rv_{B_{p}} = K_{p}\lv u_{n} \rv_{\infty}^{p-1} \rightarrow 0\quad\text{as}\quad n\rightarrow\infty. \qedhere \] 
\end{proof}

\begin{proof}[Proof of Theorem~{\normalfont{\ref{thmsmalllargeideals}\ref{itemsmall}}}]
We shall apply \Cref{thm:FSZ} with $(X ,Y) = (E,D)$, that is, either $(X,Y) = (B_p,\ell_p)$ for some $1<p<\infty$ or $(X,Y)=(S_p,c_0)$ for some $1\le p<\infty$, endowed with the $1$-UFDDs obtained by blocking the unit vector bases as follows:
\begin{equation}\label{itemsmall:eq1}  X_{n} = \operatorname{span}\{e_{j} : j \in J_{n}\}\qquad\text{and}\qquad Y_{n} = \operatorname{span}\{d_{j} : j \in J_{n}\}\qquad (n \in \N), \end{equation}
where $J_1<J_2<\cdots$ are the intervals defined in \Cref{ConstructionMPB}, and $T = \iota\in\mathscr{B}(E,D)$ is the formal inclusion map. 

Condition~\ref{thm:FSZ:i} of \Cref{thm:FSZ} is trivially satisfied because $\iota (e_j) = d_j$ for $j\in\N$.
We claim that the infimum in condition~\ref{thm:FSZ:ii} 
equals~$1$. To prove that, we begin by recalling that $Q_N = P_{L_N}$ for every $N\subseteq\N$, where~$Q_N\in\mathscr{B}(E)$ denotes the projection associated with the $1$-UFDD $(X_n)_{n\in\N}$ of~$E$, the set~$L_N$ is given by~\eqref{eq:defnLM}, and~$P_{L_N}\in\mathscr{B}(E)$ is the corresponding basis projection, as usual. Hence the claim will follow provided that we show that 
\begin{equation*}
\operatorname{dist}(\iota P_{L_M},\langle \iota P_{L_N}\rangle) = 1\qquad (M,N \in [\N],\,\lvert M \setminus N\rvert = \infty). \end{equation*}
The inequality~$\le$ is trivial because $\lVert \iota P_{L_M}\rVert = 1$. We shall verify the opposite inequality by showing that if $\operatorname{dist}(\iota P_{L_M},\langle \iota P_{L_N}\rangle) <1$ for some $M,N \in [\N]$, then $\lvert M \setminus N\rvert<\infty$. Hence, suppose  that $\lv \iota P_{L_{M}} - R\rv < 1$ for some operator $R\in\langle \iota P_{L_N}\rangle$, 
say  $R = \sum_{j=1}^{k}U_{j}\iota P_{L_{N}}V_{j}$, where $k\in\N$, $U_{1},\ldots,U_{k}\in \mathscr{B}(D)$ and $V_{1},\ldots,V_{k} \in \mathscr{B}(E)$.
By replacing $U_{j}$ with $\lv V_{j} \rv U_{j}$ and $V_{j}$ with $\frac{V_{j}}{\lv V_{j} \rv}$ if $\lv V_{j} \rv > 0$, we may suppose that $\lv V_{j} \rv \leq 1$ for each $j \in \{1,\ldots,k\}$. 
 
Take $m \in L_{M}$. Since $e_{m}$ and $d_{m} = \iota P_{L_{M}}e_m$ are unit vectors, we have
\[  \lv \iota P_{L_{M}}-R\rv\geq \lv (\iota P_{L_{M}}-R)e_{m}\rv_{D}\geq \lv \iota P_{L_{M}}e_m\rv_D -\lv Re_{m}\rv_{D} = 1 - \lv Re_{m} \rv_{D},
\]
so
\begin{equation}\label{itemsmall:eq5}
    1 - \lv \iota P_{L_{M}}-R\rv \leq \lv Re_{m} \rv_{D} \leq \sum_{j=1}^{k}\lv U_{j} \rv\, \lv \iota P_{L_{N}}V_{j}e_{m} \rv_{D} \leq k \cdot \max_{1 \leq j \leq k}\lv U_{j} \rv \cdot \lv \iota P_{L_{N}}V_{\varphi(m)}e_{m} \rv_{D}, 
\end{equation}
where we have chosen $\varphi(m) \in \{1,\ldots,k\}$ such that
\[  \max_{1 \leq j \leq k}\lv \iota P_{L_{N}}V_{j}e_{m} \rv_{D} = \lv \iota P_{L_{N}}V_{\varphi(m)}e_{m} \rv_{D}.
\]
This defines a map $\varphi\colon L_{M} \rightarrow \{1,\ldots,k\}$ which in view of~\eqref{itemsmall:eq5} satisfies
\begin{equation}\label{itemsmall:eq3}
   \lv \iota P_{L_{N}}V_{\varphi(m)}e_{m} \rv_{D}\geq \eta\qquad (m \in L_{M}),\qquad \text{where}\qquad \eta = \frac{1 - \lv \iota P_{L_{M}}-R\rv}{\displaystyle{k\cdot\max_{1 \leq j \leq k}\lv U_{j} \rv}} > 0.
\end{equation}

We use this map to introduce a new operator 
\[ W = P_{L_{N}}\sum_{j \in \varphi(L_{M})}V_{j}P_{\varphi^{-1}(\{j\})}|_{E_{L_M}}\in\mathscr{B}(E_{L_{M}},E_{L_{N}}). 
\]
Our aim is to show that it satisfies
\begin{equation}\label{itemsmall:eq4}
\inf_{m \in L_{M}}\lVert We_{m}\rVert_\infty\ge\begin{cases}
    \eta\ &\text{for}\ E=S_p,\\[1mm] {\displaystyle{\Bigl(\frac{\eta^{p}}{K_{p}}\Bigr)^{\frac{1}{p-1}}}}\ &\text{for}\ E=B_p,
\end{cases}
\end{equation}
where $K_p>0$ denotes the constant from \Cref{ThmgGJOJp}. 
Take $m \in L_{M}$, and observe that $\lVert We_{m}\rVert_\infty = \lVert P_{L_{N}}V_{\varphi(m)}e_{m}\rVert_\infty$ because 
\[ P_{\varphi^{-1}(\{j\})}e_{m} = \begin{cases}
        e_{m}\ &\text{if}\ j = \varphi(m),\\
        0\ &\text{otherwise.}
    \end{cases} \]
If $E=S_p$, then $D=c_0$, so $\lVert P_{L_{N}}V_{\varphi(m)}e_{m}\rVert_\infty = \lv \iota P_{L_{N}}V_{\varphi(m)}e_{m} \rv_{D}\ge\eta$ by~\eqref{itemsmall:eq3}, which establishes~\eqref{itemsmall:eq4} in the first case. Otherwise $E = B_{p}$ and $D=\ell_p$; combining~\eqref{itemsmall:eq3} with Jameson's inequality stated in the last line of \Cref{ThmgGJOJp}, we obtain
\begin{equation}\label{itemsmall:eq6}
\eta^p\le\lVert \iota P_{L_{N}}V_{\varphi(m)}e_{m}\rVert_{\ell_p}^p\le K_p\lVert  P_{L_{N}}V_{\varphi(m)}e_{m}\rVert_\infty^{p-1}  \lVert P_{L_{N}}V_{\varphi(m)}e_{m}\rVert_{B_p}\le K_p\lVert We_{m}\rVert_\infty^{p-1},
\end{equation}
where the simple estimate $\lVert P_{L_{N}}V_{\varphi(m)}e_{m}\rVert_{B_p}\le \lVert P_{L_{N}}\rVert\,\lVert V_{\varphi(m)}\rVert\,\lVert e_{m}\rVert_{B_p}\le 1$ justifies the final inequality. The second case of~\eqref{itemsmall:eq4} follows by rearranging~\eqref{itemsmall:eq6}.

Hence the operator~$W$ satisfies condition~\ref{GLThm1iii} of \Cref{GLThm}, so $\gl(L_{M},L_{N}) < \infty$, and therefore $\lvert M \setminus N\rvert < \infty$ by \Cref{corollaryonpageB}, as required. 

We have thus verified both conditions of \Cref{thm:FSZ}. It follows that the map~\eqref{thm:FSZ:eq1} is in\-jec\-tive. Composing it with the injection~\eqref{idealsetup}, we obtain $2^{\mathfrak{c}}$ many closed ideals of~$\mathscr{B}(E)$. They are con\-tained in the ideal of strictly singular operators because the operator~$\iota$ is strictly singular, as we showed in \Cref{P:formalinclusionmap}.
\end{proof}

\appendix
\section{Jameson's inequality for the Schreier and Baernstein norms}\label{App:GJOJ}
\noindent 
The aim of this appendix is to establish an inequality which relates the $\ell_p$-norm, the first Schreier (or $p^{\text{th}}$ Baern\-stein) norm and the $\ell_\infty$-norm. (Recall that we denote the latter by $\lVert\,\cdot\,\rVert_\infty$.) Its proof is due to Graham Jameson; we are very grateful for his permission to include it here. The inequality plays a key role in the proof of \Cref{thmsmalllargeideals}\ref{itemsmall} that we gave in \Cref{S:6}. 

\begin{theorem}[Jameson]\label{ThmgGJOJp} For every $1< p<\infty$, there is a constant $K_p\in\bigl[\frac{2^p-1}{2^{p-1}-1}, \frac{3\cdot2^{p-1}-2}{2^{p-1}-1}\bigr]$ such that 
\begin{equation}\label{Thmgraham:eq1:V2}
   \lVert x \rVert_{\ell_p}^p\leq K_p \lVert x\rVert_{\infty}^{p-1}\lVert x \rVert_{S_{1}}\qquad (x\in\K^\N).  
\end{equation}
Consequently,  $\lVert x \rVert_{\ell_p}^p\leq K_p \lVert x\rVert_{\infty}^{p-1}\lVert x \rVert_{B_{p}}$ for every $x\in B_p$.
\end{theorem}

We begin with a lemma that will help us reduce to the case of decreasing sequences.

\begin{lemma}\label{grahamlemma1V2}
Let $x\colon \N\to [0,\infty)$ be decreasing with limit~$0$. Then $\lVert x\rVert_{S_p}\le \lVert x\circ\sigma\rVert_{S_p}$ for every $1\le p<\infty$ and every permutation $\sigma\colon\N\to\N$. 
\end{lemma}

\begin{proof}
Take $F\in\mathcal{S}_1\setminus\{\emptyset\}$, and let $k=\min F$. Since~$\sigma$ is surjective, the set $\sigma^{-1}([1,2k)\cap\N)\setminus[1,k)$ contains a subset~$G$ of cardinality~$k$. Then $G\in\mathcal{S}_1$, and therefore
\[ \lVert x\circ\sigma\rVert_{S_p}^p\ge 
\sum_{n\in G} x(\sigma(n))^p\ge \sum_{j=k}^{2k-1} x(j)^p\ge\mu_p(x,F)^p, \]
where the second inequality follows because $\sigma(G)$ is a $k$-element subset of $[1,2k)\cap\N$ and~$x$ is decreasing. Now the conclusion follows by taking the supremum over~$F$.    
\end{proof}

\begin{proof}[Proof of Theorem {\normalfont{\ref{ThmgGJOJp}}}] 
Since all three norms in~\eqref{Thmgraham:eq1:V2} depend only on the moduli of the coordinates of~$x$, it suffices to consider non-negative~$x$. We may also suppose that $\lVert x\rVert_{S_1}<\infty$, as otherwise the inequality is trivial. This implies that $x\in c_0$, which in turn means that we can find a permutation $\sigma\colon\N\to\N$ such that $x\circ\sigma$ is decreasing. Therefore we can apply \Cref{grahamlemma1V2} to $x\circ\sigma$ and the permutation $\sigma^{-1}$ to obtain that $\lVert x\circ\sigma\rVert_{S_1}\le \lVert (x\circ\sigma)\circ\sigma^{-1}\rVert_{S_1} = \lVert x\rVert_{S_1}$, while 
 $\lVert x\circ\sigma \rVert_{\ell_p} =  \lVert x \rVert_{\ell_p}$ and $\lVert x\circ\sigma \rVert_{\infty} =  \lVert x \rVert_{\infty}$. In conclusion, this shows that it suffices to consider the case where $x\colon\N\to[0,\infty)$ is decreasing, and after scaling, we may suppose that  $\lVert x\rVert_{S_1}=1$.

Take $n\in\N_0$, and let $F_n = [2^n, 2^{n+1})\cap\N\in\mathcal{S}_1$. Since~$x$ is non-negative and decreasing, we have $\lVert x\rVert_\infty = x(1)$ and $ x(2^{n+1})\le x(j)\le x(2^n)$ for $j\in F_{n}$. Combining this with the fact that $\mu_1(x,F_n)\leq \lVert x\rVert_{S_1} = 1$, we obtain  
\begin{equation}\label{GJOJ:eq1}
x(2^{n+1}) \leq \frac{1}{2^{n}}\qquad\text{and}\qquad
    \mu_p(x,F_n)^p\leq x(2^n)^{p-1}\mu_1(x,F_n)\leq x(2^n)^{p-1}. \end{equation} 
Write $x(1) = \frac{\theta }{2^k}$, where $k\in\N$ and $1 \leq \theta \leq 2$.  Then we have 
\[ \sum_{n=0}^{k-1} \mu_p(x,F_n)^p = \sum_{j=1}^{2^k-1}x(j)^p\leq  (2^k-1) x(1)^p\le \theta x(1)^{p-1}. \]
Furthermore, \eqref{GJOJ:eq1} implies that $\mu_p(x,F_k)^p\leq x(1)^{p-1}$ and 
\[ \sum_{n=k+1}^\infty \mu_p(x,F_n)^p\leq  \sum_{n=k+1}^\infty \frac{1}{2^{(n-1)(p-1)}} = \Bigl(\frac{1}{2^k}\Bigr)^{p-1}\frac{1}{1-\frac{1}{2^{p-1}}} = 
\frac{x(1)^{p-1}2^{p-1}}{\theta^{p-1}(2^{p-1}-1)}. \]
Since $(F_n)_{n=0}^\infty$ is a partition of~$\N$, we conclude that 
\[ \lVert x\rVert_{\ell_p}^p = \sum_{n=0}^\infty\mu_p(x,F_n)^p\le f(\theta)x(1)^{p-1} = f(\theta)\lVert x\rVert_\infty^{p-1},\quad\text{where}\quad f(\theta ) = \theta + 1 + \frac{2^{p-1}}{(2^{p-1}-1)\theta^{p-1}}. \]   
This defines a smooth function $f\colon(0,\infty)\to(1,\infty)$ whose second derivative $f''(\theta) = \frac{2^{p-1}p(p-1)}{(2^{p-1}-1)\theta^{p+1}}$ is positive. Hence~$f$ is convex, so $\max\{f(\theta) : 1\le\theta\le 2\} = \max\{f(1),f(2)\}$. We find that \[ f(1)=f(2) = \frac{3\cdot 2^{p-1}-2}{2^{p-1}-1}, \]  
and therefore the inequality~\eqref{Thmgraham:eq1:V2} is satisfied for some constant $K_p\le\frac{3\cdot 2^{p-1}-2}{2^{p-1}-1}$.

To verify that this constant is at least~$\frac{2^p-1}{2^{p-1}-1}$, take $k\in\N$ and define $x\colon\N\to(0,\infty)$ by
\[ x(j) = \begin{cases} {\displaystyle{\frac{1}{2^k}}}\ &\text{for}\ 1 \leq j < 2^{k+1},\\[2mm] 
  {\displaystyle{\frac{1}{2^n}}}\ &\text{for}\ 2^n\le j<2^{n+1},\ \text{where}\ n\in(k,\infty)\cap\N. \end{cases} \]
Then $\lVert x\rVert_\infty = \frac{1}{2^k}$ and
\[ \lVert x\rVert_{\ell_p}^p = \frac{2^{k+1}-1}{2^{kp}} + \sum_{n=k+1}^\infty\frac{2^n}{2^{np}} = \frac{2}{2^{k(p-1)}} - \frac{1}{2^{kp}} + \frac{1}{2^{k(p-1)}(2^{p-1}-1)}. \]

We claim that $\lVert x\rVert_{S_1} = 1$. Since~$x$ is decreasing, it suffices to consider Schreier sets of the form $[j,2j)\cap\N$ for $j\in\N$ when computing~$\lVert x\rVert_{S_1}$.
Clearly $\mu_1(x,[j,2j)\cap\N)= {j}/{2^k}\le 1$ for $1\le j\le2^k$. Otherwise $j = 2^n+m$ for some $n\ge k$ and $1\le m\le 2^n$, and we have
\begin{equation*} \mu_1(x,[j,2j)\cap\N) = \frac{2^n - m}{2^n} + \frac{2m}{2^{n+1}} = 1. \end{equation*}
This proves the claim. Hence
\begin{align*} K_p\ge \frac{\lVert x\rVert_{\ell_p}^p}{\lVert x\rVert_\infty^{p-1}\lVert x\rVert_{S_1}} &= \Bigl(\frac{2}{2^{k(p-1)}} - \frac{1}{2^{kp}} + \frac{1}{2^{k(p-1)}(2^{p-1}-1)}\Bigr)2^{k(p-1)}\\ &= 2 - \frac{1}{2^k} + \frac{1}{2^{p-1}-1}
\rightarrow 2 + \frac{1}{2^{p-1}-1} = \frac{2^p-1}{2^{p-1}-1}\qquad\text{as}\qquad k\rightarrow\infty. \end{align*} 

The inequality stated in the last line of the theorem follows immediately from~\eqref{Thmgraham:eq1:V2} because the $p^{\text{th}}$ Baernstein norm $1$-dominates the first Schreier norm due to the fact that 
\[ \mu_1(x,F) = \beta_p(x,\{F\})\qquad (x\in\K^\N,\, F\in\mathcal{S}_1). \qedhere \]
\end{proof}

\end{document}